\documentclass[A4paper,10pt]{article}

\usepackage[numbers]{natbib}
\usepackage{graphicx}
\usepackage{algorithm,algorithmic}
\usepackage{setspace}
\usepackage{rotating}
\usepackage{color}
\usepackage{amsmath}
\usepackage{enumerate}
\usepackage{hyperref}
\usepackage{authblk}
\usepackage{geometry}
\usepackage{indentfirst}
\usepackage{amssymb}
\usepackage{subfigure}
\usepackage{nicefrac}
\usepackage{bm}
\usepackage{multirow}
\usepackage{amsfonts,amsthm}
\usepackage[title]{appendix}
\usepackage{epstopdf}
\usepackage{mathrsfs}
\usepackage{url}

\newcommand{\dx}{\Delta x}

\newcommand{\dt}{\Delta t}

\newcommand{\U}{\bm{U}}

\newcommand{\F}{\bm{F}}

\newcommand{\K}{\bm{K}}
\newcommand{\R}{\bm{R}}
\newcommand{\B}{\bm{B}}

\newcommand{\bH}{\bm{H}}
\newcommand{\mbH}{\mathcal{\bH}}

\newcommand{\bPsi}{\bm{\Psi}}
\newcommand{\mG}{\mathcal{G}}

\newcommand{\half}{\frac{1}{2}}
\newcommand{\nhalf}{\nicefrac{1}{2}}
\newcommand{\jmh}{{j-\frac{1}{2}}}
\newcommand{\jph}{{j+\frac{1}{2}}}

\newcommand{\dint}{\displaystyle\int}

\newcommand{\CFL}{{\rm CFL}}
\newcommand{\dvert}[1]{\vert #1 \vert}
\newcommand{\dVert}[1]{\left\Vert #1 \right\Vert}
\newcommand{\zrho}[1]{z_{#1}\rho_{#1}} 

\newtheorem{theorem}{\hspace{1mm}Theorem}[section]
\newtheorem{lemma}{\hspace{1mm}Lemma}[section]
\newtheorem{example}{\hspace{1mm}Example}[section]

\newtheorem{remark}{\hspace{1mm}Remark}[section]

\numberwithin{equation}{section}

\makeatletter
\newcommand{\doublewidetilde}[1]{{%
  \mathpalette\double@widetilde{#1}%
}}
\newcommand{\double@widetilde}[2]{%
  \sbox\z@{$\m@th#1\widetilde{#2}$}%
  \ht\z@=.9\ht\z@
  \widetilde{\box\z@}%
}

\graphicspath{{Figures/}}

\title{Adaptive mesh methods for hyperbolic conservation laws with bound-preserving flux limiters}

\author[1]{Yaguang Gu
\thanks{Email: guyaguang@scut.edu.cn}}

\author[2,3]{Guanghui Hu
\thanks{Email: garyhu@um.edu.mo}}

\author[4,5]{Tao Tang
\thanks{Email: ttang@bnbu.edu.cn}}

\affil[1]{School of Mathematics, South China University of Technology, Guangzhou 510641, China}

\affil[2]{State Key Laboratory of Internet of Things for Smart City and Department of Mathematics, University of Macau, Macao SAR 999078, China}
\affil[3]{Zhuhai UM Science $\&$ Technology Research Institute, Zhuhai 519031, China}

\affil[4] {School of Mathematics and Statistics, Guangzhou Nanfang College, Guangzhou 510970, China}

\affil[5]{Institute for Advanced Study, Beijing Normal-Hong Kong Baptist University, Zhuhai 519087, China}
\date{}

\begin{document}

\maketitle

\begin{abstract}
In this paper, we develop bound-preserving (BP) finite-volume schemes for hyperbolic conservation laws on adaptive moving meshes. For scalar conservative laws, we rewrite the conventional high-order discretization as a convex combination of first-order counterparts on each sub-cell, which is mathematically equivalent to introducing a bound-preserving flux limiter. Such a limiter is inexpensive to evaluate, with a feature that the corresponding BP CFL conditions depend solely on the first-order sub-cell schemes. A mild CFL restriction
is derived under which high-order spatial accuracy is retained.
 The proposed BP schemes are extend to two nonlinear systems, namely, the Euler equations and the five-equation transport model of two-medium flows. Numerical results demonstrate that the present schemes possess high resolution and strong robustness
properties.
\end{abstract}

\noindent
{\bf Key words:} Hyperbolic conservation laws, adaptive mesh methods, finite-volume methods, bound-preserving method.

\medskip
\noindent
{\bf AMS subject classification:} 35L65, 65M08, 65M50, 76M12.

\section{Introduction}

We are interested in developing high-order finite-volume methods  on adaptive moving meshes (AMM-FVMs) for hyperbolic conservation laws:
\begin{equation}\label{EQ:HYP_System}
u_t + f(u)_x = 0,
\end{equation}
where $u$ is a conservative variable, and $f(u)$ is the flux.
It has been demonstrated that the adaptive moving mesh methods have great potential towards enhancing the computational efficiency; see, e.g., \cite{Li2006Moving,
Tang2003Adaptive,tang2005moving}.
The main idea of the adaptive moving mesh methods is to dynamically concentrate time-dependent grid points in regions where the solution exhibits rapid variations, while coarsening the mesh in smooth areas. In this way the large gradients and discontinuities can be captured effectively. We refer the readers to some recent developments in this direction  \cite{Fu2019Arbitrary,Gu2017A,GU2022An,Kurganov2021Adaptive, Kurganov2022Well,PAN2020109558} and references therein.

Due to the interpolation errors, however, high-order AMMs usually require additional processing to ensure that the numerical solutions satisfy specific physical properties. An important property of the entropy solution of \eqref{EQ:HYP_System} is that the solution satisfies a strict maximum principle
\begin{equation*}
u_m \leq u(x,t) \leq u_M,
\end{equation*}
if $u_m \leq u(x,0) \leq u_M$ at the initial state. Therefore, the solution produced by a numerical scheme is expected to belong to the set
\begin{equation*}
\mG = \{u\,\vert\, u_m \leq u \leq u_M\},
\end{equation*}
which is called a physics-constraint admissible set, or an invariant domain. A numerical scheme satisfying the property that $u_j^{n+1}\in\mG$ whenever $u_j^n\in\mG$ is called a maximum-principle-preserving scheme, or a bound-preserving (BP) scheme. Note that the solution of hyperbolic systems in general does not admit the maximum principle. However, the concept of BP property is still crucial for numerical methods, as hyperbolicity of the governing equations needs to be maintained at the discrete level. For example, a numerical scheme is expected to preserve positive density, bounded volume fraction (between $0$ and $1$), etc. Therefore, designing a BP scheme is indispensable and has become a hot topic in the past decades.

Basically there are two BP frameworks. One is the Zhang-Shu framework \cite{ZHANG20103091}. In this framework, a convex decomposition of the high-order scheme is performed with help of numerical integration. The BP property is fulfilled in two aspects: one is to employ polynomial rescaling limiters such that reconstruction values at quadrature points are BP, and the other one is to enforce all terms in the decomposition to be BP by using provably suitable BP CFL conditions. The accuracy of the scheme is preserved as the polynomial rescaling procedure is accuracy-preserving. Later, the authors extend the idea to Euler equations \cite{ZHANG20108918,ZHANG20122245} and system with source terms \cite{ZHANG20111238}. We refer the readers to \cite{CHENG2014143,CUI2025114189,FAN2022111446,QIAN2018172,Wu2021Uniformly} and references therein for applications of this framework to a variety of hyperbolic problems.

Another type of BP framework uses the flux limiters. Assume that the conventional FVM is used for spatial discretization, that is
\begin{equation*}
\overline{u}_j^{n+1} = \overline{u}_j^n - \frac{\dt}{\dx} \left(
f_\jph - f_\jmh
\right),
\end{equation*}
Then a weighted average of fluxes computed by a low-order scheme and by a
high-order scheme is used to reduce numerical diffusions produced by the low-order scheme. The resulting numerical flux at cell boundary is given by
\begin{equation}\label{EQ:FCT}
f_\jph\approx f_\jph^{\mathrm{H}, \mathrm{(mod)}} = \theta_\jph f_\jph^{\mathrm{H}} + (1-\theta_\jph) f_\jph^{\mathrm{L}},
\end{equation}
where $f_\jph^{\mathrm{H}}$ and $f_\jph^{\mathrm{L}}$ are high-order and low-order approximations of $f_\jph$. This type of method seems to be first proposed in \cite{BOOK1975248,BORIS1976397,BORIS197338,ZALESAK1979335}, which was named the  flux-corrected transport (FCT) method. It was then further developed according to different ways in computing $\theta_\jph$; consult \cite{kuzmin2012flux} and references therein.

When the flux limiter is applied, one has to answer an intuitive question: Does the blending process in \eqref{EQ:FCT} affect the formal spatial accuracy of high-order numerical schemes? In \cite{XU2014Parametrized}, a parametrized flux limiter was proposed, where the BP property is fulfilled by using linear programming method to determine $\theta_\jmh$ and $\theta_\jph$; consult also \cite{Liang2014Parametrized,
Xiong2016Parametrized}.
It was proved in \cite{XU2014Parametrized} that the best accuracy-preserving (AP) CFL condition is $\nicefrac{1}{\sqrt{12}}$ for third-order FVMs with parametrized flux limiters.
In \cite{HU2013169}, another simpler flux limiter was proposed to enforce positivity of density and pressure in Euler equations. In this method, the whole scheme is decomposed into a convex combination
\begin{equation*}
\overline{u}_j^{n+1} = \frac{1}{2} \left( \overline{u}_j^n + 2\frac{\dt}{\dx} f_\jmh \right) + \frac{1}{2}\left( \overline{u}_j^n - 2\frac{\dt}{\dx} f_\jph \right)=:\frac{1}{2}\Phi_1 + \frac{1}{2}\Phi_2.
\end{equation*}
Then $\theta_\jph$ are determined such that $\Phi_1$ and $\Phi_2$ are positivity-preserving (PP), which is a weak version of the BP requirement. Besides, existence of AP CFL conditions can be analyzed.
Combining advantages of the methods in \cite{HU2013169,XU2014Parametrized}, a BP method was developed in \cite{Fu2025Bound,Gao2023High,GU2023112190,Gu2021A} using auxiliary function (for scalar equation) or vector (for system). For equation \eqref{EQ:HYP_System}, this method provides another convex decomposition of the scheme which reads as
\begin{equation*}
\overline{u}_j^{n+1} = \frac{1}{2} \left[ \overline{u}_j^n + 2\frac{\dt}{\dx} \left(f_\jmh - f(\overline{u}_j)\right) \right] + \frac{1}{2}\left[ \overline{u}_j^n - 2\frac{\dt}{\dx} \left(f_\jph - f(\overline{u}_j)\right) \right]=:\frac{1}{2}\widetilde{\Phi}_1 + \frac{1}{2}\widetilde{\Phi}_2,
\end{equation*}
with $f(\overline{u}_j)$ being an auxiliary function. Since $\widetilde{\Phi}_1$ and $\widetilde{\Phi}_2$ are sub-cell schemes, their BP property can be analyzed easily when the first-order numerical fluxes are employed. Similar ideas can be traced back to \cite{Guermond2017Invariant,KUZMIN2020112804} in the framework of continuous finite element methods; see also recent progress in \cite{Abgrall2025bound,Duan2025Active,
Kuzmin2023Property}.

In this paper, we will develop and analyze BP methods on adaptive moving meshes. We emphasize that the derivation is highly non-trivial as the mesh grids become time-dependent. We will focus on the FVMs which evolve the solution directly from the old mesh to the new one \cite{Fu2019Arbitrary,FU2022111600}, where the mesh redistribution is carried out based on variational approach \cite{Kurganov2021Adaptive,Tang2003Adaptive}. The paper consists three parts. Firstly, we introduce in \S\ref{Sec:BP-FVM-Scalar} a simple decomposition strategy, under which the BP property of the sub-cell schemes can be proved under certain CFL conditions. Any high-order sub-cell scheme which exceeds the required bounds can be modified by blending it with the first-order counterparts. We point out that the blending process introduces flux limiters so that the resulting finite volume schemes on adaptive moving meshes is conservative. The second part devotes to the analysis on the accuracy-preserving (AP) CFL conditions. We will prove that the AP CFL condition is $\nicefrac{1}{6}$. Thirdly, we extend the proposed BP method to two nonlinear systems: fully conservative Euler equations (in \S\ref{Sec:Euler}) and the five-equation transport model \cite{ALLAIRE2002577} of two-medium flows  which contains non-conservative product (in \S\ref{Sec:5EQ}).
Finally, numerical examples are provided in \S\ref{Sec:Num} to validate the theoretical analysis and demonstrate the robustness and effectiveness of the proposed BP AMM-FVMs.

\section{Bound-preserving method on moving meshes}\label{Sec:BP-FVM-Scalar}

In this section, we will develop an efficient BP method for \eqref{EQ:HYP_System} using finite-volume methods on moving meshes. We first introduce the conventional third-order finite-volume discretization on adaptive moving meshes. Then in \S\ref{Sec:BP-CFL} a simple bound-preserving scheme will be introduced, and its corresponding CFL conditions ensuring BP property will be obtained.  In \S\ref{Sec:AP-CFL}, accuracy-preserving CFL conditions will be studied. Adaptive mesh redistribution and WENO reconstruction will be investigated in \S\ref{Sec:mesh_redistribution} and \S\ref{Sec:WENO}, respectively. 
Let the physical domain $\Omega$ be covered with time-dependent finite volume cells $C_j(t) = [x_\jmh(t), x_\jph(t)]$ of size $\dx_j(t)$, centered at $x_j(t) = (x_\jmh(t)+x_\jph(t))/2$. Movements of grid points are time-dependent, but we assume that each grid point moves along a straight line with a uniform grid velocity $\omega_\jph$ during each time step, namely
\begin{equation*}
x_\jph(t) = x_\jph^{n}+\omega_\jph (t-t^n), \quad t\in[t^n, t^{n+1}].
\end{equation*}
The computation of the grid velocity will be introduced in \S\ref{Sec:mesh_redistribution}. To avoid impractically small control volume, we also assume that there exists a constant $C$ independent of $\dx_j$ such that $\dx\leq C\dx_j$ for all $j$, where $\dx=\max_j\{\dx_j\}$.

Assume that the solution, realized in terms of its cell average $\overline{u}_j$ over $C_j$, is available at a certain time level $t^n$. We integrate the governing equation \eqref{EQ:HYP_System} over $[x_\jmh(t), x_\jph(t)]$, and use the Reynolds transport theorem and the divergence theorem to obtain
\begin{equation}\label{EQ:ODE-u}
\begin{aligned}
\frac{\mathrm{d}}{\mathrm{d}t} \left[ \dx_j(t)\overline{u}_j(t) \right] =  \frac{\mathrm{d}}{\mathrm{d}t} \int_{x_\jmh(t)}^{x_\jph(t)} u\;\mathrm{d}x &= \int_{x_\jmh(t)}^{x_\jph(t)} \left[ u_t + (\omega u)_x \right]\;\mathrm{d}x \\
&= -\int_{x_\jmh(t)}^{x_\jph(t)} (f(u) -\omega u)_x \;\mathrm{d}x = -\left[ H_\jph(t) - H_\jmh(t) \right],
\end{aligned}
\end{equation}
where $H_\jph(t):= H(\omega_\jph,u_\jph(t))= f(u_\jph(t))-\omega_\jph u_\jph(t)$.
The above ODE is discretized by the conventional finite-volume method in space and the third-order strong stability-preserving Runge-Kutta (SSP-RK3) method \cite{GOTTLIEB1998Total,Shu2001Strong} in time:
\begin{subequations}\label{EQ:u-RK}
\begin{align}
\overline{u}_j^{(1)} &=
\frac{\dx_j^n}{\dx_j^{(1)}}\overline{u}_j^n - \frac{\dt^n}{\dx_j^{(1)}}\left[
\mbH \left(\omega_\jph,u_\jph^{n,-},u_\jph^{n,+}\right) - \mbH \left(\omega_\jmh,u_\jmh^{n,-},u_\jmh^{n,+}\right)\right],\label{EQ:u-1st-Stage-RK} \vspace{+3pt}\\
\overline{u}_j^{(2)} &= \frac{3}{4}\frac{\dx_j^n}{\dx_j^{(2)}}\overline{u}_j^n + \frac{1}{4}\left\lbrace
\frac{\dx_j^{(1)}}{\dx_j^{(2)}}\overline{u}_j^{(1)} - \frac{\dt^n}{\dx_j^{(2)}}\left[
\mbH \left(\omega_\jph,u_\jph^{(1),-},u_\jph^{(1),+}\right) - \mbH \left(\omega_\jmh,u_\jmh^{(1),-},u_\jmh^{(1),+}\right)
\right]\right\rbrace,\label{EQ:u-2nd-Stage-RK} \vspace{+3pt}\\
\overline{u}_j^{n+1} &= \frac{1}{3}\frac{\dx_j^n}{\dx_j^{n+1}}\overline{u}_j^n + \frac{2}{3}\left\lbrace
\frac{\dx_j^{(2)}}{\dx_j^{n+1}}\overline{u}_j^{(2)} - \frac{\dt^n}{\dx_j^{n+1}}\left[
\mbH \left(\omega_\jph,u_\jph^{(2),-},u_\jph^{(2),+}\right) - \mbH \left(\omega_\jmh,u_\jmh^{(2),-},u_\jmh^{(2),+}\right)
\right]\right\rbrace,\label{EQ:u-3rd-Stage-RK}
\end{align}
\end{subequations}
where $u_\jph^{(\ell),\pm}$, $\ell = 0, 1, 2$ (with $u_\jph^{(0),\pm}=u_\jph^{n,\pm}$), are obtained by the third-order WENO reconstruction, and $\mbH$ is any Lipschitz continuous, monotone, and consistent numerical flux. For illustration purpose, we take the simple Lax-Friedrichs (LF) flux as follows,
\begin{equation}\label{LF1}
\mbH( \omega, u^-, u^+ ) = \frac{1}{2}(H(\omega,u^-) + H(\omega,u^+)) - \frac{\alpha}{2}(u^+ - u^-),
\end{equation}
where $\alpha = \max\dvert{f'(u)-\omega}$. The time step size is restricted by the CFL condition, that is
\begin{equation}\label{EQ:Conventional-dt}
\dt^n = \CFL\frac{ \min\left\lbrace\dx_j^n\right\rbrace}{\alpha^n}.
\end{equation}
Conventionally, $\CFL\leq 1$ is taken when LF flux is applied.

While evolving the solution from $t^n$ to $t^{n+1}$, intermediate mesh sizes also need to be updated. Meanwhile, the discrete geometric conservation laws (D-GCL), which states that a uniform flow does not change with time, should be exactly preserved at the discrete level in the presence of arbitrary mesh movement. 
To this end, we assume that $u=u^0=\mathrm{const}$. 
Then the ODE \eqref{EQ:ODE-u} becomes
\begin{equation}\label{EQ:ODE-dx}
\frac{\mathrm{d}}{\mathrm{d}t} \left[ \dx_j(t)u^0 \right] = (w_\jph - w_\jmh) u^0.
\end{equation}
Eliminating the constant $u^0$ from \eqref{EQ:ODE-dx}, and employing the same SSP-RK3 method for the resulting ODE, we obtain the intermediate mesh sizes used in \eqref{EQ:u-RK}.
\begin{equation}\label{EQ:dx-RK}
\begin{aligned}
\dx_j^{(1)} &=
\dx_j^n + \dt^n (\omega_\jph - \omega_\jmh), \vspace{+3pt}\\
\dx_j^{(2)} &=
\frac{3}{4}\dx_j^n + \frac{1}{4}\dx_j^{(1)} + \frac{1}{4}
\dt^n (\omega_\jph - \omega_\jmh), \vspace{+3pt}\\
\dx_j^{n+1} &=
\frac{1}{3}\dx_j^n + \frac{2}{3}\dx_j^{(2)} + \frac{2}{3}
\dt^n (\omega_\jph - \omega_\jmh).
\end{aligned}
\end{equation}
By using \eqref{EQ:u-RK} and \eqref{EQ:dx-RK} the D-GCL is satisfied, which will be numerically verified in example \ref{exam-GCL} in \S\ref{Sec:Num}.

\subsection{Bound-preserving method and its CFL conditions}
\label{Sec:BP-CFL}

Due to reconstruction errors, the numerical scheme (\ref{EQ:u-RK})-(\ref{LF1}) may not satisfy the maximum principle.
In this subsection, we will design a bound-preserving scheme by employing flux limiters, and then derive the corresponding BP CFL conditions. The following lemma is useful for deriving the BP method.

\begin{lemma}\label{Lem:HmH}
Let $\omega$ be the grid velocity, and let $u$ be the conservative variable. Then we have
\begin{equation}\label{EQ:HmH}
\mbH(\omega,u,u) - \mbH(0,u,u) = -\omega u.
\end{equation}
\end{lemma}
\begin{proof}
The consistency of the numerical flux reads as $\mbH(\omega,u,u) = f(u)-\omega u$, by which we have
\begin{equation*}
\mbH(\omega,u,u) - \mbH(0,u,u) = \left(f(u)-\omega u\right) - f(u) = -\omega u.
\end{equation*}
This finishes the proof of Lemma \ref{Lem:HmH}.
\end{proof}

To illustrate the idea, let us first consider the forward Euler temporal discretization, which coincides with the first stage of the SSP-RK3 method. Using Lemma \ref{Lem:HmH}, we can reformulate \eqref{EQ:u-1st-Stage-RK} to obtain
\begin{equation*}
\begin{aligned}
\overline{u}_j^{(1)} =&
\frac{\dx_j^n}{\dx_j^{(1)}}\overline{u}_j^n - \frac{\dt}{\dx_j^{(1)}}\left[
\mbH \left(\omega_\jph,u_\jph^{n,-},u_\jph^{n,+}\right) - \mbH \left(\omega_\jmh,u_\jmh^{n,-},u_\jmh^{n,+}\right)\right] \vspace{+3pt}\\
=& \frac{\dx_j^n}{\dx_j^{(1)}}\overline{u}_j^n - \frac{\dt}{\dx_j^{(1)}}\left[
\mbH \left(\omega_\jph,u_\jph^{n,-},u_\jph^{n,+}\right)
- \mbH \left(\omega_\jph,\overline{u}_j^n,\overline{u}_j^n\right)
+ \mbH \left(\omega_\jph,\overline{u}_j^n,\overline{u}_j^n\right)
- \mbH \left(0,\overline{u}_j^n,\overline{u}_j^n\right) \right. \vspace{+3pt}\\
&\left.
- \mbH \left(\omega_\jmh,u_\jmh^{n,-},u_\jmh^{n,+}\right)
+ \mbH \left(\omega_\jmh,\overline{u}_j^n,\overline{u}_j^n\right)
- \mbH \left(\omega_\jmh,\overline{u}_j^n,\overline{u}_j^n\right)
+ \mbH \left(0,\overline{u}_j^n,\overline{u}_j^n\right)\right] \vspace{+3pt}\\
\overset{\eqref{EQ:HmH}}{=}& \overline{u}_j^n - \frac{\dt^n }{\dx_j^{(1)}}\left\lbrace\left[
\mbH \left(\omega_\jph,u_\jph^{n,-},u_\jph^{n,+}\right) - \mbH \left(\omega_\jph,\overline{u}_j^n,\overline{u}_j^n\right) \right] - \left[\mbH \left(\omega_\jmh,u_\jmh^{n,-},u_\jmh^{n,+}\right) - \mbH \left(\omega_\jmh,\overline{u}_j^n,\overline{u}_j^n\right) \right]\right\rbrace.
\end{aligned}
\end{equation*}
Then we decompose $\overline{u}_j^{(1)}$ into a convex combination of sub-cell schemes
\begin{equation}\label{EQ:Convex-Decomp-RK1}
\overline{u}_j^{(1)} = \frac{1}{2} u_j^{(1),\rm H,-} + \frac{1}{2}u_j^{(1),\rm H,+},
\end{equation}
where
\begin{eqnarray}
&& u_j^{(1),\rm H,-} = \overline{u}_j^n + 2\lambda_j^n\left[\mbH
(\omega_\jmh,u_\jmh^{n,-},u_\jmh^{n,+}) - \mbH(\omega_\jmh,\overline{u}_j^n,\overline{u}_j^n)\right], \label{2x1}\\
&& u_j^{(1),\rm H,+} = \overline{u}_j^n - 2\lambda_j^n\left[\mbH
(\omega_\jph,u_\jph^{n,-},u_\jph^{n,+}) - \mbH(\omega_\jph,\overline{u}_j^n,\overline{u}_j^n)\right],\label{2x2}
\end{eqnarray}
with
\begin{equation*}
\lambda_j^n = \frac{\dt^n}{\dx_j^{(1)}} = \frac{\dt^n}{\dx_j^n + \dt^n(\omega_\jph-\omega_\jmh)}.
\end{equation*}
By convexity of $\mG$, it follows from (\ref{EQ:Convex-Decomp-RK1}) that
\[
 \overline{u}_j^{(1)}\in\mG \quad {\rm if} \quad u_j^{(1),\rm H,\pm}\in\mG.
 \]
Unfortunately, for standard CFL condition, $u_j^{(1),\rm H,-}$ and $u_j^{(1),\rm H,+}$ may not belong to $\mG$. Therefore, a refined CFL condition has to be considered. To this end , we consider the first-order counterparts of (\ref{2x1})-(\ref{2x2}):
\begin{equation*}
\begin{aligned}
u_j^{(1),\rm L,-} &= \overline{u}_j^n + 2\lambda_j^n\left[\mbH
(\omega_\jmh,\overline{u}_{j-1}^n,\overline{u}_j^n) - \mbH(\omega_\jmh,\overline{u}_j^n,\overline{u}_j^n)\right], \\
u_j^{(1),\rm L,+} &= \overline{u}_j^n - 2\lambda_j^n\left[\mbH
(\omega_\jph,\overline{u}_j^n,\overline{u}_{j+1}^n) - \mbH(\omega_\jph,\overline{u}_j^n,\overline{u}_j^n)\right].
\end{aligned}
\end{equation*}
We show in the following lemma that both $u_j^{(1),\rm L,-}$ and $u_j^{(1),\rm L,+}$ are BP under certain CFL conditions.

\begin{lemma}\label{Lem:BP-CFL-RK1}
Let $\overline{u}_j^n\in\mG$ be given for all $j$. Then $u_j^{(1),\rm L,+}\in\mG$, provided that
\begin{equation}\label{EQ:BP-CFL-RK1}
\lambda_j^{n}\alpha^{n} = \frac{\dt^n \alpha^{n}}{\dx_j^n + \dt^n(\omega_\jph-\omega_\jmh)} \leq \frac{1}{2}.
\end{equation}
where $\alpha^n$ is the LF parameter at the first stage.
Similarly, $u_j^{(1),\rm L,-}\in\mG$ under the same CFL condition.
\end{lemma}

\begin{proof}
We first derive the BP CFL condition for $u_j^{(1),\rm L,+}\in\mG$.
Substituting the Lax-Friedrichs numerical flux into the expression of $u_j^{(1),\rm L,+}$, we obtain that
\begin{equation} \label{2x3}
u_j^{(1),\rm L,+} = \Phi_1(\overline{u}_j^n) + \Phi_2(\overline{u}_{j+1}^n),
\end{equation}
where
\begin{equation} \label{2x4}
\Phi_1(u) = \left[ 1- \lambda_j^n(\alpha^n+\omega_\jph)\right] u + \lambda_j^n f(u), \qquad
\Phi_2(u) = \lambda_j^n(\alpha^n+\omega_\jph) u - \lambda_j^n f(u).
\end{equation}
Inspired by \cite{ZHANG20103091}, we seek for conditions under which $u_j^{(1),\rm L,+}$ is monotone increasing with respect to both $\overline{u}_{j}^n$ and $\overline{u}_{j+1}^n$. By examining (\ref{2x3})-(\ref{2x4}), this requires that
\begin{eqnarray*}
&& \Phi_1'(\overline{u}_j^n) = 1- \lambda_j^n\left(\alpha^n+\omega_\jph - f'(\overline{u}_{j}^n)\right) \geq 0,  \\
&& \Phi_2'(\overline{u}_{j+1}^n) = \lambda_j^n\left(\alpha^n+\omega_\jph - f'(\overline{u}_{j+1}^n)\right) \geq 0.
\end{eqnarray*}
The second inequality holds immediately due to the definition of the parameter $\alpha^n$. Observe that
\begin{equation}
\dVert{\alpha^{n} - \left(f'(\overline{u}_{j}^{n}) - \omega_\jph\right)} \leq 2\alpha^{n}.
\end{equation}
Then the first inequality holds if \eqref{EQ:BP-CFL-RK1} is satisfied.
Finally, using monotonicity of $u_j^{(1),\rm L,+}$ gives
\begin{equation*}
u_m=\Phi_1(u_m)+\Phi_2(u_m)\leq u_j^{(1),\rm L,+}\leq \Phi_1(u_M)+\Phi_2(u_M)=u_M.
\end{equation*}
Hence, we have proved that $u_m\leq u_j^{(\ell),\rm L,+}\leq u_M$, i.e. $u_j^{(\ell),\rm L,+}\in\mG$. The BP CFL condition \eqref{EQ:BP-CFL-RK1} for $u_j^{(1),\rm L,-}\in\mG$ can be derived similarly, and is thus omitted for brevity. This completes the proof of the lemma.
\end{proof}

With the help of $u_j^{(1),\rm L,\pm}\in\mG$, we modify $u_j^{(1),\rm H,\pm}$ as follows. We blend $u_j^{(1),\rm H,\pm}$ with their first-order counterparts to obtain
\begin{equation}\label{EQ:blend-RK1}
\widetilde{u}_{j}^{(1),\rm H,-} = \theta_\jmh u_{j}^{(1),\rm H,-} + (1-\theta_\jmh)u_{j}^{(1),\rm L,-}, \qquad
\widetilde{u}_j^{(1),\rm H,+} = \theta_\jph u_j^{(1),\rm H,+} + (1-\theta_\jph)u_j^{(1),\rm L,+},
\end{equation}
and then replace $u_j^{(1),\rm H,\pm}$ in \eqref{EQ:Convex-Decomp-RK1} by $\widetilde{u}_j^{(1),\rm H,\pm}$. This gives
\begin{equation*}
\overline{u}_j^{(1)} = \frac{1}{2} \widetilde{u}_j^{(1),\rm H,-} + \frac{1}{2}\widetilde{u}_j^{(1),\rm H,+},
\end{equation*}
In (\ref{EQ:blend-RK1}), $\theta_\jph = \min\{ \theta_\jph^-, \theta_\jph^+ \}$, with
\begin{equation}\label{EQ:theta}
\begin{aligned}
\theta_\jph^- &= \min\left\lbrace 1,
~\Theta(u_j^{(1),\rm H,+}-u_m, u_j^{(1),\rm L,+}-u_m), ~\Theta(u_M-u_j^{(1),\rm H,+}, u_M-u_j^{(1),\rm L,+})\right\rbrace,  \vspace{+3pt}\\
\theta_\jph^+ &= \min\left\lbrace 1,
~\Theta(u_{j+1}^{(1),\rm H,+}-u_m, u_{j+1}^{(1),\rm L,+}-u_m), ~\Theta(u_M-u_{j+1}^{(1),\rm H,+}, u_M-u_{j+1}^{(1),\rm L,+})
\right\rbrace,
\end{aligned}
\end{equation}
where 
\begin{equation}\label{EQ:Operator_Theta}
\Theta(\phi^{\rm H}, \phi^{\rm L}) = \bigg|\frac{\phi^{\rm L}-\varepsilon^\mathrm{BP}}{\phi^{\rm L}-\phi^{\rm H}}\bigg|,
\end{equation}
with the small positive number $\varepsilon^\mathrm{BP}$ being introduced to avoid rounding errors. In numerical examples, we take $\varepsilon^\mathrm{BP}=10^{-16}$.
The parameters $\theta_\jph$ ensure that $\widetilde{u}_{j}^{(1),\rm H,\pm}\in\mG$ for all $j$.

We remark that we have introduced a common parameter $\theta_\jph$ in \eqref{EQ:blend-RK1} at the cell boundary. This is due to the following observations.
We reformulate $\widetilde{u}_j^{(1),\rm H,+}$ and $\widetilde{u}_{j+1}^{(1),\rm H,-}$ to obtain
\begin{eqnarray*}
&& \widetilde{u}_j^{(1),\rm H,+} = \overline{u}_j^n - 2\lambda_j^n\left[\widetilde{\mbH}_\jph^n - \mbH(\omega_\jph,\overline{u}_j^n,\overline{u}_j^n)\right], \\
&& \widetilde{u}_{j+1}^{(1),\rm H,-} = \overline{u}_{j+1}^n + 2\lambda_{j+1}^n\left[\widetilde{\mbH}_\jph^n - \mbH(\omega_\jph,\overline{u}_{j+1}^n,\overline{u}_{j+1}^n)\right],
\end{eqnarray*}
where
\begin{equation*}
\widetilde{\mbH}_\jph^n := \theta_\jph\mbH
(\omega_\jph,u_\jph^{n,-},u_\jph^{n,+}) + (1-\theta_\jph)\mbH
(\omega_\jph,\overline{u}_{j}^n,\overline{u}_{j+1}^n).
\end{equation*}
It reveals that a common modified flux $\widetilde{\mbH}_\jph^n$ is employed at the cell boundary. Therefore, the scheme is still conservative after enforcing BP property. In addition, the blending process introduces in essence flux limiters, that is, the modified flux $\widetilde{\mbH}_\jph^n$ is obtained by blending the high-order flux with its first-order counterpart, with the blending parameter $\theta_\jph$ being introduced in \eqref{EQ:theta} to enforce BP property.

Next, we consider the $\ell$-th stage of the SSP-RK3 method, which can be written as
\begin{eqnarray}
&& \dx_j^{(\ell)}\overline{u}_j^{(\ell)} = \xi^{(\ell)}\dx_j^n\overline{u}_j^n \label{EQ:RK-stage-l} \\
&& \qquad  \quad + (1-\xi^{(\ell)})\left\lbrace
\dx_j^{(\ell-1)}\overline{u}_j^{(\ell-1)} - \dt^n\left[
\mbH \left(\omega_\jph,u_\jph^{(\ell-1),-},u_\jph^{(\ell-1),+}\right) - \mbH \left(\omega_\jmh,u_\jmh^{(\ell-1),-},u_\jmh^{(\ell-1),+}\right)
\right]\right\rbrace, \nonumber
\end{eqnarray}
where $\dx_j^{(\ell)}$ is computed by
\begin{equation}\label{EQ:RK-stage-dx-l}
\dx_j^{(\ell)} = \xi^{(\ell)}\dx_j^n + (1-\xi^{(\ell)})\left[
\dx_j^{(\ell-1)} + \dt(\omega_\jph - \omega_\jmh)\right].
\end{equation}
Here, $\ell = 1, 2, 3$, with $\overline{u}_j^{(0)}=\overline{u}_j^{n}$, $\overline{u}_j^{(3)}=\overline{u}_j^{n+1}$, $\xi^{(1)} = 0$,
$\xi^{(2)} = \frac{3}{4}$,
$\xi^{(3)} = \frac{1}{3}$.
Using the same trick as before, we rewrite \eqref{EQ:RK-stage-l} as
\begin{equation*}
\begin{array}{rl}
\overline{u}_j^{(\ell)} =& \dfrac{\xi^{(\ell)}\dx_j^n}{\dx_j^{(\ell)}}\overline{u}_j^n + \dfrac{(1-\xi^{(\ell)})}{\dx_j^{(\ell)}}\left\lbrace
\dx_j^{(\ell-1)}\overline{u}_j^{(\ell-1)} - \dt^n\left[
\left(\mbH
(\omega_\jph,u_\jph^{(\ell-1),-},u_\jph^{(\ell-1),+}) - \mbH(\omega_\jph,\overline{u}_j^{(\ell-1)},\overline{u}_j^{(\ell-1)})\right)\right.\right.  \vspace{+3pt}\\
&+\left.\left. \left(\mbH(\omega_\jph,\overline{U}_j^{(\ell-1)},\overline{u}_j^{(\ell-1)}) - \mbH
(0,\overline{u}_j^{(\ell-1)},\overline{u}_j^{(\ell-1)})\right) - \left(
\mbH
(\omega_\jmh,u_\jmh^{(\ell-1),-},u_\jmh^{(\ell-1),+}) - \mbH(\omega_\jmh,\overline{u}_j^{(\ell-1)},\overline{u}_j^{(\ell-1)})
\right)\right.\right. \vspace{+3pt}\\
&-\left.\left.\left(
\mbH
(\omega_\jmh,u_\jmh^{(\ell-1),-},u_\jmh^{(\ell-1),+}) - \mbH(0,\overline{u}_j^{(\ell-1)},\overline{u}_j^{(\ell-1)})
\right)
\right]\right\rbrace.
\end{array}
\end{equation*}
Then we use Lemma \ref{Lem:HmH} and \eqref{EQ:RK-stage-dx-l} to deliver a decomposition of \eqref{EQ:RK-stage-l}:
\begin{equation}\label{EQ:Convex-Decomp}
\overline{u}_j^{(\ell)} = \dfrac{\xi^{(\ell)}\dx_j^n}{\dx_j^{(\ell)}}\overline{u}_j^n + \dfrac{(1-\xi^{(\ell)})}{\dx_j^{(\ell)}}\left[
\dx_j^{(\ell-1)} + \dt^n(\omega_\jph - \omega_\jmh)\right]\left(
\frac{1}{2} u_j^{(\ell),\rm H,-} + \frac{1}{2} u_j^{(\ell),\rm H,+}
\right),
\end{equation}
where the sub-cell schemes $u_j^{(\ell),\rm H,\pm}$ read as
\begin{equation}\label{EQ:ujH-l}
\begin{array}{rl}
u_j^{(\ell),\rm H,-} = \overline{u}_j^{(\ell-1)} + 2\lambda_j^{(\ell-1)}\left[\mbH
(\omega_\jmh,u_\jmh^{(\ell-1),-},u_\jmh^{(\ell-1),+}) - \mbH(\omega_\jmh,\overline{u}_j^{(\ell-1)},\overline{u}_j^{(\ell-1)})\right], \vspace{+3pt}\\
u_j^{(\ell),\rm H,+} = \overline{u}_j^{(\ell-1)} - 2\lambda_j^{(\ell-1)}\left[\mbH
(\omega_\jph,u_\jph^{(\ell-1),-},u_\jph^{(\ell-1),+}) - \mbH(\omega_\jph,\overline{u}_j^{(\ell-1)},\overline{u}_j^{(\ell-1)})\right].
\end{array}
\end{equation}
with
\begin{equation}\label{EQ:Lambda-ell-1}
\lambda_j^{(\ell-1)} = \frac{\dt^n}{\dx_j^{(\ell-1)}+\dt^n(\omega_\jph - \omega_\jmh)}.
\end{equation}
Then we modify $u_j^{(\ell),\rm H,\pm}$ utilizing their first-order counterparts $u_j^{(\ell),\rm L,\pm}$ to obtain
\begin{equation}\label{EQ:blend-RKl}
\widetilde{u}_{j}^{(\ell),\rm H,-} = \theta_\jmh u_{j}^{(\ell),\rm H,-} + (1-\theta_\jmh)u_{j}^{(\ell),\rm L,-}, \qquad
\widetilde{u}_j^{(\ell),\rm H,+} = \theta_\jph u_j^{(\ell),\rm H,+} + (1-\theta_\jph)u_j^{(\ell),\rm L,+},
\end{equation}
where $u_j^{(\ell),\rm L,\pm}$ read as
\begin{equation}\label{EQ:ujL-l}
\begin{array}{rl}
u_j^{(\ell),\rm L,-} = \overline{u}_j^{(\ell-1)} + 2\lambda_j^{(\ell-1)}\left[\mbH
(\omega_\jmh,\overline{u}_{j-1}^{(\ell-1)},\overline{u}_{j}^{(\ell-1)}) - \mbH(\omega_\jmh,\overline{u}_j^{(\ell-1)},\overline{u}_j^{(\ell-1)})\right], \vspace{+3pt}\\
u_j^{(\ell),\rm L,+} = \overline{u}_j^{(\ell-1)} - 2\lambda_j^{(\ell-1)}\left[\mbH
(\omega_\jph,\overline{u}_{j}^{(\ell-1)},\overline{u}_{j+1}^{(\ell-1)}) - \mbH(\omega_\jph,\overline{u}_j^{(\ell-1)},\overline{u}_j^{(\ell-1)})\right].
\end{array}
\end{equation}
The conditions for $u_j^{(\ell),\rm L,\pm}\in\mG$ will be given later in Theorem \ref{Th:BP-Cond}. Here, the blending parameter
\[
\theta_\jph = \min\{ \theta_\jph^-, \theta_\jph^+ \},
 \]
with $\theta_\jph^\pm$ being defined the same way as in \eqref{EQ:theta} and \eqref{EQ:Operator_Theta}.
Similar to the forward Euler method, we can replace $u_{j}^{(\ell),\rm H,\pm}$ by $\widetilde{u}_{j}^{(\ell),\rm H,\pm}$ in \eqref{EQ:Convex-Decomp}, or modify the numerical flux at $x_\jph$ by
\begin{equation}\label{EQ:flux-RKl-blending}
\widetilde{\mbH}_\jph^{(\ell-1)} := \theta_\jph\mbH
(\omega_\jph,u_\jph^{(\ell-1),-},u_\jph^{(\ell-1),+}) + (1-\theta_\jph)\mbH
(\omega_\jph,\overline{u}_{j}^{(\ell-1)},\overline{u}_{j+1}^{(\ell-1)}).
\end{equation}

Taking all stages of SSP-RK3 into account, we give BP CFL conditions in the following theorem.

\begin{theorem}\label{Th:BP-Cond}
Let $\overline{u}_j^n\in\mG$ be given for all $j$, where $\overline{u}_j^n$ is the numerical solution of SSP-RK3 described above. Then $\overline{u}_j^{n+1}\in\mG$ for all $j$ if the following conditions are satisfied:
\begin{enumerate}[(1)]
\item The intermediate lengths $\dx_j^{(\ell-1)}$ (with $\ell=1,2,3$, $\dx_j^{(0)} = \dx_j^n$, $\dx_j^{(2)} = \dx_j^{n+1}$) satisfy that, for any $j$ and $\ell$,
\begin{equation}\label{EQ:dxj-Ineq}
\dx_j^{(\ell-1)}+\dt^n(\omega_\jph - \omega_\jmh) \geq 0;
\end{equation}
\item The time step size $\dt^n$ satisfies that, for any $j$ and $\ell$,
\begin{equation}\label{EQ:BP-CFL}
\lambda_j^{(\ell-1)}
\alpha^{(\ell-1)} = \frac{\dt^n \alpha^{(\ell-1)}}{\dx_j^{(\ell-1)}+\dt^n(\omega_\jph - \omega_\jmh)} \leq \frac{1}{2}.
\end{equation}
\end{enumerate}
\end{theorem}
\begin{proof}
The first stage of SSP-RK3 ensures that $\overline{u}_j^{(1)}\in\mG$, which has been proved in Lemma \ref{Lem:BP-CFL-RK1}. In the following, we prove that $\overline{u}_j^{(\ell)}\in\mG$ for $\ell=2,3$ by reduction. Suppose that $\overline{u}_j^{(\ell-1)}\in\mG$.
If \eqref{EQ:BP-CFL} is satisfied, then we can prove that $u_j^{(\ell),\rm L,\pm}\in\mG$ using the same trick as that used in Lemma \ref{Lem:BP-CFL-RK1}. In addition, note that if the condition \eqref{EQ:dxj-Ineq} is satisfied then $\overline{u}_j^{(\ell)}$ in \eqref{EQ:Convex-Decomp} is a convex combination of $\overline{u}_j^n$, $u_j^{(\ell),\rm L,-}$ and $u_j^{(\ell),\rm L,+}$. As a result, we have $\overline{u}_j^{(\ell)}\in\mG$ by convexity of the admissible set $\mG$. This completes the proof of the theorem.
\end{proof}

\begin{remark} {\em
On uniform mesh, i.e., when $\omega_\jph=0$ for all $j$, $u_j^{(\ell),\rm H,\pm}$ and $u_j^{(\ell),\rm L,\pm}$ become
\begin{equation*}
\begin{array}{ll}
u_j^{(1),\rm H,-} = \overline{u}_j^n + 2\lambda_j^n\left[\mbH
(0,u_\jmh^{n,-},u_\jmh^{n,+}) - \mbH(0,\overline{u}_j^n,\overline{u}_j^n)\right], \quad
&u_j^{(1),\rm H,+} = \overline{u}_j^n - 2\lambda_j^n\left[\mbH
(0,u_\jph^{n,-},u_\jph^{n,+}) - \mbH(0,\overline{u}_j^n,\overline{u}_j^n)\right], \vspace{+3pt}\\
u_j^{(1),\rm L,-} = \overline{u}_j^n + 2\lambda_j^n\left[\mbH
(0,\overline{u}_{j-1}^n,\overline{u}_j^n) - \mbH(0,\overline{u}_j^n,\overline{u}_j^n)\right], \quad
&u_j^{(1),\rm L,+} = \overline{u}_j^n - 2\lambda_j^n\left[\mbH
(0,\overline{u}_j^n,\overline{u}_{j+1}^n) - \mbH(0,\overline{u}_j^n,\overline{u}_j^n)\right],
\end{array}
\end{equation*}
where $\mbH(0,u_1,u_2)$ is the classical LF flux, namely
\begin{equation*}
\mbH(0,u_1,u_2) = \frac{1}{2}(f(u_1) + f(u_2)) - \frac{\alpha}{2}(u_2-u_1), \quad
\alpha = \max_u\{\dvert{f'(u)}\}.
\end{equation*}
The corresponding BP CFL condition becomes
\begin{equation*}
\lambda_j^{(\ell-1)}
\alpha^{(\ell-1)} = \alpha^{(\ell-1)}\frac{\dt^n}{\dx} \leq \frac{1}{2}.
\end{equation*}
In addition, it is worth pointing out that the BP method on uniform mesh (i.e., $\omega_\jph=0$) coincides with the one proposed in our previous studies \cite{Fu2025Bound,Gao2023High,GU2023112190,Gu2021A}; see also monolithic convex limiting method in \cite{Kuzmin2023Property,WONG2021A}.
}
\end{remark}

\subsection{Accuracy-preserving CFL conditions}\label{Sec:AP-CFL}
The bound-preserving method proposed in the previous section is carried out by blending the high-order numerical flux with its first-order counterpart. The aim of this subsection is to answer an intuitive question: Does the blending process affect the spatial accuracy of the scheme? In another word, does
\begin{equation}\label{EQ:Err-H-m-H}
\dVert{\widetilde{\mbH}_\jph^{(\ell-1)} - \mbH(\omega_\jph,u_\jph^{(\ell-1),-},u_\jph^{(\ell-1),+)} } = \mathcal{O}(\dx^3)
\end{equation}
hold under reasonable accuracy-preserving (AP) CFL conditions?
In this section, we will answer this question by deriving the accuracy-preserving CFL conditions.

To simplify the discussion, we only consider the flux limiter for lower bound, as the flux limiter for upper bound can be analyzed similarly. We substitute $\theta_\jph$ into \eqref{EQ:flux-RKl-blending}, and make use of the definition of $\Theta$ defined in \eqref{EQ:Operator_Theta} to obtain
\begin{equation*}
\begin{array}{l}
\dVert{\widetilde{\mbH}_\jph^{(\ell-1)} - \mbH(\omega_\jph,u_\jph^{(\ell-1),-},u_\jph^{(\ell-1),+})} \vspace{+3pt}\\
= \begin{cases}
\frac{\dVert{\overline{u}_j-2\lambda_j^{(\ell-1)}\left[
\mbH(\omega_\jph,u_\jph^{(\ell-1),-},u_\jph^{(\ell-1),+}) - \mbH(\omega_\jph,\overline{u}_j,\overline{u}_j)
 \right]-u_m}}{2\lambda_j^{(\ell-1)}}, &\theta_\jph = \Theta(u_j^{(\ell),\rm H,+}-u_m, u_j^{(\ell),\rm L,+}-u_m), \vspace{+3pt}\\
\frac{\dVert{\overline{u}_{j+1}+2\lambda_{j+1}^{(\ell-1)}\left[
\mbH(\omega_\jph,u_\jph^{(\ell-1),-},u_\jph^{(\ell-1),+}) - \mbH(\omega_\jph,\overline{u}_{j+1},\overline{u}_{j+1})
\right]-u_m}}{2\lambda_{j+1}^{(\ell-1)}}, &\theta_\jph = \Theta(u_{j+1}^{(\ell),\rm H,-}-u_m, u_{j+1}^{(\ell),\rm L,-}-u_m).
\end{cases}
\end{array}
\end{equation*}
Here, we have excluded the discussion on $\theta_\jph = 1$, since the BP limiter is not activated in this situation. Moreover, we ignore the influence of the parameter $\varepsilon^\mathrm{BP}$ in our analysis.

We first derive AP CFL conditions when $u_{j}^{(\ell),\rm H,+}<u_m\leq u_{j}^{(\ell),\rm L,+}$. In this situation, we have
\begin{equation*}
\theta_\jph = \frac{u_{j}^{(\ell),\rm L,+}-u_m}{u_{j}^{(\ell),\rm L,+}-u_{j}^{(\ell),\rm H,+}}.
\end{equation*}
Since at any stage of SSP-RK3,
\begin{equation*}
u_\jph^{\pm} = u(x_\jph) + \mathcal{O}(\dx^3), \qquad
\overline{u}_j = \frac{1}{\dx_j}\int_{C_j} u(x)\,\mathrm{d}x + \mathcal{O}(\dx^3),
\end{equation*}
where $u(x)$ is the exact solution,
it suffices to analyze at the continuous level. Notice that the quantities above are time-dependent, but we temporally omit this dependence for the sake of brevity. To obtain the main result Theorem \ref{Th:AP-CFL}, we need the following lemmas.

\begin{lemma}\label{Lem:AP-CFL-R}
Let $u(x)$ be the solution of (\ref{EQ:HYP_System}). Assume that $\frac{\partial H}{\partial u}(\omega,u) = f'(u) - \omega_\jph$ is bounded, $u_{j}^{(\ell),\rm H,+}<u_m\leq u_{j}^{(\ell),\rm L,+}$, and
\begin{equation}\label{EQ:Scheme-less-um-R}
\overline{u}_j - 2\lambda_j\left(
\mbH(\omega_\jph,u_\jph,u_\jph) - \mbH(\omega_\jph,\overline{u}_j,\overline{u}_j)\right) < u_m.
\end{equation}
If the CFL condition
\begin{equation}\label{EQ:CFL-Cond-R}
\lambda_j \max_{x\in C_j}\dVert{f'(u(x)) - \omega_\jph} \leq \frac{1}{6}
\end{equation}
is satisfied, then we have
\begin{equation}\label{EQ:Err-HmH-R}
\frac{1}{2\lambda_j}\dVert{\overline{u}_j - 2\lambda_j\left[
\mbH(\omega_\jph,u_\jph,u_\jph) - \mbH(\omega_\jph,\overline{u}_j,\overline{u}_j)\right] -u_m} = \mathcal{O}(\dx^3).
\end{equation}
\end{lemma}
\begin{proof}
By consistency of the flux, we have
\begin{equation*}
\overline{u}_j - 2\lambda_j\left(
\mbH(\omega_\jph,u_\jph,u_\jph) - \mbH(\omega_\jph,\overline{u}_j,\overline{u}_j)\right)
=\overline{u}_j - 2\lambda_j\left[
\left(f(u_\jph)-\omega_\jph u_\jph\right) - \left(f(\overline{u}_j)-\omega_\jph \overline{u}_j\right)
\right].
\end{equation*}
In the following, we will prove case by case that
\begin{equation*}
\overline{u}_j - 2\lambda_j\left[
\left(f(u_\jph)-\omega_\jph u_\jph\right) - \left(f(\overline{u}_j)-\omega_\jph \overline{u}_j\right)\right]\geq u_m + \mathcal{O}(\dx_j^3),
\end{equation*}
which indicates, together with \eqref{EQ:Scheme-less-um-R}, that
\begin{equation}\label{EQ:Err-Scheme-um-R}
\dVert{\overline{u}_j - 2\lambda_j\left[
\left(f(u_\jph)-\omega_\jph u_\jph\right) - \left(f(\overline{u}_j)-\omega_\jph \overline{u}_j\right)\right]-u_m} = \mathcal{O}(\dx_j^3).
\end{equation}

{\bf Case I}: there exists a global minimum point $x_m$ belonging to closure of $C_j$, or a local minimum point $x_m$ inside $C_j$, s.t. $u'(x_m)=0$, and $u''(x_m)>0$. Under this assumption, we have $\dvert{x-x_m}=\mathcal{O}(\dx_j)$, $\forall x\in C_j$. We perform Taylor expansion of $u(x)$ at $x_m$ to obtain
\begin{equation*}
u(x) = u(x_m) + \frac{1}{2}u''(x_m)(x-x_m)^2 + \mathcal{O}(\dx_j^3).
\end{equation*}
Then we integrate $u(x)$ over $C_j$ to get estimate of the cell average over $C_j$
\begin{equation}\label{EQ:u-bar}
\overline{u}_j = u(x_m) + \frac{1}{2} u''(x_m)\left((x_j-x_m)^2+\frac{1}{12}\dx_j^2\right) + \mathcal{O}(\dx_j^3),
\end{equation}
and evaluate $u(x)$ at $x_\jph$ to obtain value at cell interface
\begin{equation*}
u_\jph = u(x_m) + \frac{1}{2} u''(x_m)\left(
(x_j-x_m)^2 + (x_j-x_m)\dx_j + \frac{1}{4}\dx_j^2
\right) + \mathcal{O}(\dx_j^3).
\end{equation*}
Afterwards, we substitute $\overline{u}_j$ and $u_\jph$ into $f(u)$ to find that
\begin{equation*}
f(\overline{u}_j) = f(u(u_m)) + \frac{1}{2}u''(x_m)f'(u(x_m))\left((x_j-x_m)^2+\frac{1}{12}\dx_j^2\right) + \mathcal{O}(\dx_j^3),
\end{equation*}
and
\begin{equation}\label{EQ:f-u-jph}
f(u_\jph) = f(u(u_m)) + \frac{1}{2}u''(x_m)f'(u(x_m))\left(
(x_j-x_m)^2 + (x_j-x_m)\dx_j + \frac{1}{4}\dx_j^2
\right) + \mathcal{O}(\dx_j^3).
\end{equation}
Combining \eqref{EQ:u-bar}-\eqref{EQ:f-u-jph} yields
\begin{equation*}
\overline{u}_j - 2\lambda_j\left[
\left(f(u_\jph)-\omega_\jph u_\jph\right) - \left(f(\overline{u}_j)-\omega_\jph \overline{u}_j\right)\right] = u(x_m) + \frac{1}{2}u''(x_m)R(x_m) + \mathcal{O}(\dx_j^3),
\end{equation*}
where the function $R$ is a quadratic function of $x_m$, and takes the form
\begin{equation*}
R(x_m) = (x_j-x_m)^2+\frac{1}{12}\dx_j^2 -2\lambda_j\left(f'(u(x_m))-\omega_\jph\right)\left(
(x_j-x_m)\dx_j + \frac{1}{6}\dx_j^2
\right).
\end{equation*}
The function $R$ reaches its minimum at $x_m^c = x_j - \lambda_j(f'(u(x_m))-\omega_\jph)\dx_j$. Subsequently, we compute minimum value of $R$:
\begin{equation*}
R(x_m^c) = -\frac{\dx_j^2}{12}\left[6\lambda_j\left(f'(u(x_m))-\omega_\jph\right)-1\right]\left[2\lambda_j\left(f'(u(x_m))-\omega_\jph\right)+1\right].
\end{equation*}
Under the CFL condition \eqref{EQ:CFL-Cond-R}, $R(x_m^c)\geq 0$. Therefore we have, with $u''(x_m)>0$ and $u(x_m)\geq u_m$, that
\begin{equation*}
\overline{u}_j - 2\lambda_j\left[
\left(f(u_\jph)-\omega_\jph u_\jph\right) - \left(f(\overline{u}_j)-\omega_\jph \overline{u}_j\right)\right] \geq u_m + \mathcal{O}(\dx_j^3).
\end{equation*}

{\bf Case II:} $u(x)$ is monotone increasing over $C_j$. In this case, $u'(x)>0$, and the local minimum point is $x_m = x_{\jmh}$. We perform Taylor expansion of $u(x)$ at $x_\jmh$ to obtain
\begin{equation*}
u(x) = u(x_\jmh) + u'(x_\jmh)(x-x_\jmh) + \frac{1}{2}u''(x_\jmh)(x-x_\jmh)^2 + \mathcal{O}(\dx_j^3).
\end{equation*}
Then we follow the same lines in Case I to get
\begin{equation*}
\begin{aligned}
u_\jph &= u_\jmh + u_\jmh'\dx_j + \frac{1}{2} u_\jmh''\dx_j^2 + \mathcal{O}(\dx_j^3), \vspace{+3pt}\\
f(u_\jph) &= f(u_\jmh) + f'(u_\jmh)u_\jmh'\dx_j + \frac{1}{2}f'(u_\jmh)u_\jmh''\dx_j^2 + \frac{1}{2}f''(u_\jmh)(u_\jmh')^2\dx_j^2 + \mathcal{O}(\dx_j^3), \vspace{+3pt}\\
\overline{u}_j &= u_\jmh + \frac{1}{2}u_\jmh'\dx_j + \frac{1}{6} u_\jmh''\dx_j^2 + \mathcal{O}(\dx_j^3), \vspace{+3pt}\\
f(\overline{u}_j) &= f(u_\jmh) + \frac{1}{2}f'(u_\jmh)u_\jmh'\dx_j + \frac{1}{6}f'(u_\jmh)u_\jmh''\dx_j^2 + \frac{1}{8}f''(u_\jmh)(u_\jmh')^2\dx_j^2 + \mathcal{O}(\dx_j^3).
\end{aligned}
\end{equation*}
Equipped with these Taylor expansions, we derive that
\begin{equation*}
\begin{array}{l}
\overline{u}_j - 2\lambda_j\left[
\left(f(u_\jph)-\omega_\jph u_\jph\right) - \left(f(\overline{u}_j)-\omega_\jph \overline{u}_j\right)\right] \vspace{+3pt}\\
= u_\jmh + r_1 u_\jmh'\dx_j + r_2 (u_\jmh')^2\dx_j^2 + \frac{1}{2}r_3u_\jmh''\dx_j^2 + \mathcal{O}(\dx_j^3), \vspace{+3pt}\\
=u(x_\jmh-\sqrt{r_3}\dx_j+r_2u'(x_\jmh-\sqrt{r_3}\dx_j)\dx_j^2) + (r_1 + \sqrt{r_3})u_\jmh'\dx_j,
\end{array}
\end{equation*}
where we have defined $\alpha_j = f'(u_\jmh)-\omega_\jph$, and
\begin{equation*}
r_1 = \frac{1}{2} - \lambda_j\alpha_j, \quad
r_2 = -\frac{3}{4}\lambda_j f''(u_\jmh), \quad
r_3 = \frac{4}{3}\left(\frac{1}{4}-\lambda_j\alpha_j\right).
\end{equation*}
Under the CFL condition \eqref{EQ:CFL-Cond-R}, $r_1\geq 0$ and $r_3\geq 0$. Hence, we have, together with $u_\jmh'>0$, that
\begin{equation*}
\overline{u}_j - 2\lambda_j\left[
\left(f(u_\jph)-\omega_\jph u_\jph\right) - \left(f(\overline{u}_j)-\omega_\jph \overline{u}_j\right)\right] \geq u_m + \mathcal{O}(\dx_j^3).
\end{equation*}

{\bf Case III:} $u(x)$ is monotone decreasing over $C_j$. In this case, $u'(x)<0$, and the local minimum point is $x_m = x_\jph$. We perform Taylor expansion of $u(x)$ at $x_\jph$ to find that
\begin{equation*}
\begin{aligned}
\overline{u}_j &= u_\jph - \frac{1}{2}u_\jph'\dx_j + \frac{1}{6}u_\jph''\dx_j^2 + \mathcal{O}(\dx_j^3), \vspace{+3pt}\\
f(\overline{u}_j) &= f(u_\jph) - \frac{1}{2}f'(u_\jph)u_\jph'\dx_j + \frac{1}{6}f'(u_\jph)u_\jph''\dx_j^2 + \frac{1}{8}f''(u_\jph)(u_\jph')^2\dx_j^2 + \mathcal{O}(\dx_j^3).
\end{aligned}
\end{equation*}
Then we derive that
\begin{equation*}
\begin{array}{l}
\overline{u}_j - 2\lambda_j\left[
\left(f(u_\jph)-\omega_\jph u_\jph\right) - \left(f(\overline{u}_j)-\omega_\jph \overline{u}_j\right)\right] \vspace{+3pt}\\
= u_\jph - r_1 u_\jph'\dx_j + r_2 (u_\jph')^2\dx_j^2 + \frac{1}{2}r_3u_\jph''\dx_j^2 + \mathcal{O}(\dx_j^3), \vspace{+3pt}\\
=u(x_\jph+\sqrt{r_3}\dx_j+r_2u'(x_\jph+\sqrt{r_3}\dx_j)\dx_j^2) - (r_1 + \sqrt{r_3})u_\jph'\dx_j,
\end{array}
\end{equation*}
where we have defined $\alpha_j = f'(u_\jph)-\omega_\jph$, and
\begin{equation*}
r_1 = \frac{1}{2} + \lambda_j\alpha_j, \quad
r_2 = \frac{1}{4}\lambda_j f''(u_\jph), \quad
r_3 = \frac{2}{3}\left(\frac{1}{2}+\lambda_j\alpha_j\right).
\end{equation*}
Under the CFL condition \eqref{EQ:CFL-Cond-R}, $r_1\geq 0$ and $r_3\geq 0$. Hence, we have, together with $u_\jmh'<0$, that
\begin{equation*}
\overline{u}_j - 2\lambda_j\left[
\left(f(u_\jph)-\omega_\jph u_\jph\right) - \left(f(\overline{u}_j)-\omega_\jph \overline{u}_j\right)\right] \geq u_m + \mathcal{O}(\dx_j^3).
\end{equation*}
Combining the above three cases, we have proved \eqref{EQ:Err-Scheme-um-R}.
Finally, \eqref{EQ:Err-HmH-R} holds since $\dx\leq C\dx_j$, $f'(u)-\omega$ is bounded, and the CFL number $\nicefrac{1}{6}$ on RHS of \eqref{EQ:CFL-Cond-R} is independent of $\dx_j$.
\end{proof}

In the case that $u_{j+1}^{(\ell),\rm H,-}<u_m\leq u_{j+1}^{(\ell),\rm L,-}$,
we employ the same way as in Lemma \ref{Lem:AP-CFL-R} to obtain the AP CFL condition below.

\begin{lemma}\label{Lem:AP-CFL-L}
Let $u(x)$ be the solution of (\ref{EQ:HYP_System}). Assume that $\frac{\partial H}{\partial u}(\omega,u) = f'(u) - \omega_\jph$ is bounded, $u_{j+1}^{(\ell),\rm H,-}<u_m\leq u_{j+1}^{(\ell),\rm L,-}$, and
\begin{equation*}
\overline{u}_{j+1} + 2\lambda_{j+1}\left(
\mbH(\omega_\jph,u_\jph,u_\jph) - \mbH(\omega_\jph,\overline{u}_{j+1},\overline{u}_{j+1})\right) < u_m.
\end{equation*}
If the CFL condition
\begin{equation*}
\lambda_{j+1} \max_{x\in C_{j+1}}\dVert{f'(u(x)) - \omega_\jph} \leq \frac{1}{6}
\end{equation*}
is satisfied, then we have
\begin{equation*}
\frac{1}{2\lambda_{j+1}}\dVert{\overline{u}_{j+1} + 2\lambda_{j+1}\left[
\mbH(\omega_\jph,u_\jph,u_\jph) - \mbH(\omega_\jph,\overline{u}_{j+1},\overline{u}_{j+1})\right] -u_m} = \mathcal{O}(\dx^3).
\end{equation*}
\end{lemma}

Using the definition of the parameter $\alpha$ in LF flux, and the definition of $\lambda_j$ defined in \eqref{EQ:Lambda-ell-1}, we conclude our AP CFL conditions below by applying the above two lemmas.

\begin{theorem}\label{Th:AP-CFL}
Assume that $\overline{u}_j^n\in\mG$ for all $j$, where $\overline{u}_j^n$ is the numerical solution of SSP-RK3 described above. If
\begin{equation}\label{EQ:AP-CFL}
\lambda_j^{(\ell-1)}\alpha^{(\ell-1)} = \frac{\dt^n \alpha^{(\ell-1)}}{\dx_j^{(\ell-1)}+\dt^n(\omega_\jph - \omega_\jmh)}
 \leq \frac{1}{6}
\end{equation}
holds for all $j$ and $\ell$, then
\eqref{EQ:Err-H-m-H} holds for all stages of SSP-RP3 at all cell interfaces, i.e., the spatial accuracy is preserved.
\end{theorem}

Note that on uniform mesh, the AP CFL condition \eqref{EQ:AP-CFL} becomes
\begin{equation*}
\lambda_j^{(\ell-1)}\alpha^{(\ell-1)} = \alpha^{(\ell-1)}\frac{\dt^n}{\dx}
 \leq \frac{1}{6}.
\end{equation*}

\subsection{Adaptive moving mesh redistribution}\label{Sec:mesh_redistribution}

Let us first assume that the mesh grids $a = x_{\nhalf}^n < x_{\nicefrac{3}{2}}^n < x_{N+\nhalf}^n = b$ of the physical domain $\Omega=[a,b]$ and the cell averages $\overline{u}_j^n$ are available at time level $t^n$.
In addition, we also introduce a uniform partition of logical domain $[0,1]$, which has same number of grid points as the physical domain, namely $0 = \xi_{\nhalf} < \xi_{\nicefrac{3}{2}} < \dots,<\xi_{N+\nhalf}=1$, with $\Delta\xi = 1/N$. A one-to-one mapping from the logical domain to the physical one
\begin{equation*}
x = x(\xi), \quad
\xi\in [0,1], \quad
x(0) = a, \quad
x(1) = b,
\end{equation*}
is determined by solving the following equation \cite{Tang2003Adaptive}:
\begin{equation}\label{EQ:Euler-Lagrange}
(\sigma x_\xi)_\xi = 0,
\end{equation}
where $\sigma$ is a monitor function depending on the differentiations of physical solution $u$, which will be specified in \S\ref{Sec:Num}. Following \cite{Kurganov2021Adaptive}, we start with $\varphi_j^0 = D\overline{u}_j^n$, and then smooth out $\varphi_j^0$ to avoid very singular mesh, by
\begin{equation*}
\varphi_j^{(\ell+1)} = \frac{1}{4}(\varphi_{j-1}^{(\ell)} + 2\varphi_j^{(\ell)} + \varphi_{j+1}^{(\ell)}), \quad
l = 0, 1, \dots, m-1.
\end{equation*}
Then the monitor function is defined by
\begin{equation}\label{EQ:sigma}
\sigma_j = 1+\alpha\varphi_j^{(m)},
\end{equation}
where the intensity parameter $\alpha$ controls concentration of mesh grids, and is computed in the following way
\begin{equation*}
\alpha = \left(\dfrac{1-\beta}{\beta(b-a)}\dint_a^b \varphi^{(m)} \,\mathrm{d}x \right)^{-1},
\end{equation*}
with $\beta\in (0,1)$ being a parameter depending on the problem at hand.

After obtaining the monitor function, we solve \eqref{EQ:Euler-Lagrange} by Jacobi iteration
\begin{equation*}
x_{j+\half}^{[\nu+1]} = \dfrac{\sigma_j x_{j-\half}^{[\nu]} + \sigma_{j+1} x_{j+\frac{3}{2}}^{[\nu]}}{\sigma_j + \sigma_{j+1}},\quad
1\leq j\leq N-1, \quad
\nu = 0,1,\dots,s-1.
\end{equation*}
where we have set $x_{j+\half}^{[0]} = x_{j+\half}^n$, and $x_{\half}^{[\nu+1]} = a$, $x_{N+\half}^{[\nu+1]} = b$. In all numerical tests, we set $m=s=8$.

Finally, we limit immoderate movement of the grid points by the following two steps.

\noindent\textbf{Step 1} (avoid twisted mesh):
For each inner grid point, we set
\begin{equation*}
\widetilde{x}_{j+\half}^{[s]} =
\begin{cases}
x_j^n, & x_{j+\half}^{[s]} < x_j^n, \vspace{+3pt}\\
x_{j+1}^n, & x_{j+\half}^{[s]} > x_{j+1}^n, \vspace{+3pt}\\
x_{j+\half}^{[s]}, & \text{ otherwise.}
\end{cases}
\end{equation*}
\noindent\textbf{Step 2} (avoid very small mesh size): If $\Delta \widetilde{x}_j^{[s]} < \dx_{\min}$, where $\dx_{\min}$ is the minimal allowed mesh size, then we modify $\widetilde{x}_{j+\half}^{[s]} = x_{j\pm\half}^{n}$, and set $\varphi_{j-1}^{[0]} = \varphi_{j}^{[0]} = \varphi_{j+1}^{[0]} = 0$ as the initial monitor function (before smoothing) at the next time step. For all test problems, we set $\dx_{\min} = \dx_\mathrm{unif}/20$, where $\dx_\mathrm{unif} = (b-a)/N$ is the uniform mesh size using same number of grid points.

After modifications, the grid velocity in time period $[t^n, t^{n+1}]$ is computed by $\omega_{j+\half} = (\widetilde{x}_{j+\half}^{[s]} - x_{j+\half}^{n})/{\dt^n}$.

\subsection{WENO reconstruction}\label{Sec:WENO}
Assume that the solution, realized in terms of its cell averages $\overline{u}_j$, are available at a certain time level. Then a WENO polynomial is reconstructed on $C_j$ following the steps below.
\begin{enumerate}
\item First, define the reconstruction patch $S_j = \{C_{j-2}, \cdots, C_{j+2}\}$.
\item Then, a quadratic polynomial $p_j^{(k)}(x)$ is reconstructed on sub-stencil $S_j^{(k)} = \{C_{j-k}, C_{j-k+1}, C_{j-k+2}\}$, $k = 0, 1, 2$, such that
\begin{equation*}
\frac{1}{\dVert{C_{\ell}}}\int_{C_{\ell}} p_j^{(k)}(x)\,\mathrm{d}x = \overline{u}_{\ell}, \quad
\ell = j-k,~ j-k+1, ~j-k+2.
\end{equation*}
\item Finally, combine polynomials $p_j^{(k)}(x)$ to obtain the eventual non-oscillatory quadratic reconstruction polynomial $p_j(x)$ by
\begin{equation*}
\operatorname{W}_j[u] : = p_j(x) = \sum_k \omega_k p_j^{(k)}(x), \quad
\omega_k = \frac{\alpha_k}{\sum_{k=0}^2 \alpha_k}, \quad
\alpha_k = \frac{d_k}{(\beta_k + \mu^2\varepsilon)^2}.
\end{equation*}
Here, $\beta_k$ is the smoothness indicator defined in \cite{JIANG1996202}:
\begin{equation*}
\beta_k = \sum_{\ell=1}^{2} \int_{C_j} (\dx_j)^{2\ell - 1} \left( \frac{\partial^\ell p_j^{(k)}(x)}{\partial x^\ell} \right)^2 \,\mathrm{d}x.
\end{equation*}
$\varepsilon = 10^{-12}$ is a small number to avoid division by zero. The weights $d_k$ of the candidate polynomial $p_j^{(k)}(x)$ are set as $d_0 = d_2 = 0.25$ and $d_1 = 0.5$.
The parameter $\mu$ takes the form
\begin{equation*}
\mu = 10^{-40} + \frac{1}{5}\sum_{\ell=j-2}^{j+2} \dVert{\overline{u}_{\ell} - \overline{u}_{S_j}}, \quad
\overline{u}_{S_j} = \frac{1}{5} \sum_{\ell=j-2}^{j+2} \overline{u}_{\ell},
\end{equation*}
which is affine-invariant, that is $\operatorname{W}[\lambda u+c] = \lambda\operatorname{W}[u] + c$, for any scale constant $\lambda$ and translation constant $c$. The affine-invariant WENO operator is capable of dealing with structures at different scales, especially at very small scale; consult \cite{Fu2025Bound,Gao2023High,WANG2022630}.
\end{enumerate}

We remark that for nonlinear systems,
the WENO reconstruction is employed in characteristic space in order to compress spurious oscillations.

To close this section, we make some remarks on some implementation issues.
In AP CFL condition \eqref{EQ:AP-CFL}, the time step size $\dt^n$ depends implicitly on the solution and on the lengths of cell values at all stages in the process of SSP-RK3.
In practice, we start with the conventional $\dt^n$ as computed in \eqref{EQ:Conventional-dt}, and follow \S\ref{Sec:mesh_redistribution} to compute grid velocity $\omega_\jph$, which will be frozen from $t^n$ to $t^{n+1}$. If the AP CFL condition is not satisfied for some $1\leq j\leq N_x$ at a certain stage $\ell$ of SSP-RK3, then we repeat the time step with a twice smaller $\dt^n$. We emphasize that every time we repeat the time step, we need to adjust the location of the grid points $x_\jph(t^{n+1})$ by the grid velocity $\omega_\jph$ and the new time step size $\dt^n$. The AP CFL condition \eqref{EQ:AP-CFL} excludes the possibility of infinite loops, since it holds as $\dt^n\to 0$. The flowchart is given in Algorithm \ref{Alg:flowchart-AP}.

\begin{algorithm}[h]
\renewcommand{\algorithmicrequire}{\textbf{INPUT:}}
\renewcommand{\algorithmicensure}{\textbf{OUTPUT:}}
\caption{Flowchart of AMM-FVM}\label{Alg:flowchart-AP}
\begin{algorithmic}[1]
\REQUIRE $x_\jph^n$, $\dx_j^n$, and BP solution $\overline{u}_j^n$,
\ENSURE $x_\jph^{n+1}$, $\dx_j^{n+1}$, and BP solution $\overline{u}_j^{n+1}$.

\STATE Compute $\dt^n$ by \eqref{EQ:Conventional-dt}.

\STATE Compute grid velocity $\omega_\jph$ following \S\ref{Sec:mesh_redistribution}. Freeze $\omega_\jph$ in the computation below.

\FOR {$\ell = 1, 2, 3$,}

\IF {the AP CFL condition \eqref{EQ:AP-CFL} is satisfied,}

\STATE Perform WENO reconstruction described in
\S\ref{Sec:WENO} to obtain $u_\jph^{(\ell-1),\pm}$.

\STATE Compute high-order schemes $\overline{u}_j^{(\ell),\mathrm{H},\pm}$ by \eqref{EQ:ujH-l} as well as their first-order counterparts $\overline{u}_j^{(\ell),\mathrm{L},\pm}$ by \eqref{EQ:ujL-l}.

\STATE Compute BP limiter $\theta_\jph$ and modify numerical flux in \eqref{EQ:flux-RKl-blending}.

\STATE Evolve solution to obtain $\overline{u}_j^{(\ell)}$ by \eqref{EQ:u-RK} using modified numerical flux, and update $\dx_j^{(\ell)}$ by \eqref{EQ:dx-RK}.

\STATE Update mesh grids by
\begin{equation*}
x_\jph^{(\ell)} = \begin{cases}
x_\jph^n + \dt^n\omega_\jph, & \ell = 1, 3,\\
x_\jph^n + \frac{\dt^n}{2}\omega_\jph, & \ell = 2,
\end{cases}
\end{equation*}
where $x_\jph^{(3)} = x_\jph^{n+1}$, which is the location of the mesh grids at time level $t^{n+1}$.

\ELSE

\STATE Set $\dt^n = \dt^n/2$. Go back to Step 4.

\ENDIF

\ENDFOR
\end{algorithmic}
\end{algorithm}

\section{Application to Euler equations}\label{Sec:Euler}
In this section, we extend our third-order BP-FVM scheme to the Euler equations of gas dynamics, which can be written in a vector form as
\begin{equation*}
\begin{array}{l}
\U_t + \F(\U)_x = 0,\qquad
\U = \begin{bmatrix}
\rho\\
\rho u\\
E
\end{bmatrix},\qquad
\F(\U) = \begin{bmatrix}
\rho u\\
\rho u^2 + p\\
(E+p)u
\end{bmatrix},
\end{array}
\end{equation*}
where $\rho$, $u$, $E$ and $p$ are the density, velocity, total energy and  pressure, respectively. The system is closed by the ideal equation of state (EOS), which describes the thermodynamic properties of a fluid, and takes the form in the case of ideal gas:
\begin{equation*}
p = (\gamma-1)\left(E-\dfrac{1}{2}\rho u^2\right),
\end{equation*}
where $\gamma$ is the specific heat ratio, taken as $1.4$ in this paper.

Similar to \eqref{EQ:u-1st-Stage-RK}, \eqref{EQ:u-2nd-Stage-RK}, and \eqref{EQ:u-3rd-Stage-RK}, the third-order finite-volume discretization on adaptive moving meshes can be rewritten as
\begin{eqnarray}
&& \overline{\U}_j^{(\ell)} = \xi^{(\ell)}\frac{\dx_j^n}{\dx_j^{(\ell)}}\overline{\U}_j^n \label{EQ:RK-stage-l-Euler}\\
&& \qquad \quad + (1-\xi^{(\ell)})\left\lbrace
\frac{\dx_j^{(\ell-1)}}{\dx_j^{(\ell)}}\overline{\U}_j^{(\ell-1)} - \frac{\dt^n}{\dx_j^{(\ell)}}\left[
\mbH \left(\omega_\jph,\U_\jph^{(\ell-1),-},\U_\jph^{(\ell-1),+}\right)-\mbH \left(\omega_\jmh,\U_\jmh^{(\ell-1),-},\U_\jmh^{(\ell-1),+}\right)
\right]\right\rbrace, \nonumber
\end{eqnarray}
where $\mbH$ is the Lax-Friedrichs numerical flux:
\begin{equation} \label{3x1}
\mbH\left(\omega_\jph,\U_\jph^{(\ell-1),-},\U_\jph^{(\ell-1),+}\right) = \frac{1}{2}\left(\bH(\U_\jph^{(\ell-1),-}) + \bH(\U_\jph^{(\ell-1),+})\right) - \frac{\alpha^{(\ell-1)}}{2}\left(\U_\jph^{(\ell-1),+} - \U_\jph^{(\ell-1),-}\right),
\end{equation}
Here, $\bH(\omega,\U) = \F(\U) - \omega\U$, and the parameter $\alpha^{(\ell-1)}$ is the spectral radius of $\frac{\partial \F}{\partial \U} - \omega \bm{I}$, which is computed by
\begin{equation*}
\alpha^{(\ell-1)} = \max_j\left\lbrace\dVert{u_j^{(\ell-1)}-\omega_\jmh} + c_j^{(\ell-1)},\; \dVert{u_j^{(\ell-1)}-\omega_\jph} + c_j^{(\ell-1)}\right\rbrace,
\end{equation*}
where $c = \sqrt{\gamma p/\rho}$ is the sound speed.
For the Euler equations, only lower bounds need to be preserved. The convex physics-constraint admissible set is
\begin{equation*}
\mG = \{\U \,\vert\, \rho(\U)>0, \, p(\U)>0 \}.
\end{equation*}
We use the same idea as in \S\ref{Sec:BP-FVM-Scalar} to rewrite the scheme \eqref{EQ:RK-stage-l-Euler} into \eqref{EQ:Convex-Decomp}, and define $\U_j^{(\ell),\rm H,\pm}$ and $\U_j^{(\ell),\rm L,\pm}$ the same way as in \eqref{EQ:ujH-l} and \eqref{EQ:ujL-l}, respectively.

In the following, we prove that $\U_j^{(\ell),\rm L,\pm}\in\mG$ under certain BP CFL conditions.

\begin{lemma}(\cite{FU2022111600})\label{Lem:BP-CFL-Euler}
If $\overline{\U}_j^{(\ell-1)}\in\mG$ for all $j$, where $\overline{\U}_j^{(\ell-1)}$ is defined by
(\ref{EQ:RK-stage-l-Euler})-(\ref{3x1}), then
\begin{equation*}
\Psi_1(\omega_{j\pm\frac{1}{2}}, \overline{\U}_j^{(\ell-1)})\in\mG, \qquad
\Psi_2(\omega_{j\pm\frac{1}{2}}, \overline{\U}_j^{(\ell-1)})\in\mG,
\end{equation*}
where
\begin{equation*}
\Psi_1(\omega, \U) := \U - \dfrac{1}{\alpha} \left( \F(\U) - \omega\U \right),  \qquad
\Psi_2(\omega, \U) := \U + \dfrac{1}{\alpha} \left( \F(\U) - \omega\U \right).
\end{equation*}
\end{lemma}

Using the above lemma, we can establish the following BP CFL condition for the Euler equations.

\begin{theorem}
Assume that $\overline{\U}_j^{(\ell-1)}\in\mG$ for all $j$, where
$\overline{\U}_j^{(\ell-1)}$ is defined by
(\ref{EQ:RK-stage-l-Euler})-(\ref{3x1}). Then $\U_j^{(\ell),\rm L,\pm}\in\mG$ under the BP CFL condition
\begin{equation}\label{EQ:BP-CFL-Euler}
\lambda_j^{(\ell-1)}\alpha^{(\ell-1)} = \frac{\dt^n \alpha^{(\ell-1)}}{\dx_j^{(\ell-1)}+\dt^n(\omega_\jph - \omega_\jmh)}
 \leq \frac{1}{2}.
\end{equation}
\end{theorem}
\begin{proof}
According to the proof of Lemma \ref{Lem:BP-CFL-RK1}, we obtain
\begin{equation*}
\begin{aligned}
\U_j^{(\ell),\rm L,-} &=
\left( 1- 2\lambda_j^{(\ell-1)}\alpha^{(\ell-1)}\right) \overline{\U}_j^{(\ell-1)} + \lambda_j^{(\ell-1)}\alpha^{(\ell-1)}\Psi_1(\omega_\jmh, \overline{\U}_j^{(\ell-1)}) + \lambda_j^{(\ell-1)}\alpha^{(\ell-1)}\Psi_2(\omega_\jmh, \overline{\U}_{j-1}^{(\ell-1)}), \vspace{+3pt}\\
\U_j^{(\ell),\rm L,+} &=
\left( 1- 2\lambda_j^{(\ell-1)}\alpha^{(\ell-1)}\right) \overline{\U}_j^{(\ell-1)} + \lambda_j^{(\ell-1)}\alpha^{(\ell-1)}\Psi_2(\omega_\jph, \overline{\U}_j^{(\ell-1)}) + \lambda_j^{(\ell-1)}\alpha^{(\ell-1)}\Psi_1(\omega_\jph, \overline{\U}_{j+1}^{(\ell-1)}).
\end{aligned}
\end{equation*}
By Lemma \ref{Lem:BP-CFL-Euler}, we have $\U_j^{(\ell),\rm L,\pm}\in\mG$, provided that they are convex combinations of BP quantities (i.e., $\overline{\U}_j^{(\ell-1)}$, $\Psi_1$, and $\Psi_2$). Note that in both $\U_j^{(\ell),\rm L,-}$ and $\U_j^{(\ell),\rm L,+}$, the sum of coefficients is $1$.
Therefore, $\U_j^{(\ell),\rm L,\pm}\in\mG$ if the coefficients are non-negative, i.e.,
\[
1- 2\lambda_j^{(\ell-1)}\alpha^{(\ell-1)}\geq 0, \quad \lambda_j^{(\ell-1)}\alpha^{(\ell-1)}\geq 0,
\]
which is true under the assumption (\ref{EQ:BP-CFL-Euler}). This completes the proof of the theorem.
\end{proof}

In practice, we first blend $\U_j^{(\ell),\rm H,\pm}$ with their first-order counterparts to obtain $\widetilde{\U}_j^{(\ell),\rm H,\pm}$, whose densities are positive. Then flux limiters are applied to $\widetilde{\U}_j^{(\ell),\rm H,\pm}$ to obtain $\doublewidetilde{\U}_j^{(\ell),\rm H,\pm}$ to ensure the positivity of the  densities and pressures.


\section{Application to five-equation transport model of two-medium flows}\label{Sec:5EQ}

In this section, we extend our BP method to hyperbolic system with non-conservative product, which can be written as
\begin{equation}\label{EQ:Hyp-NC-prod}
\U_t + \F(\U)_x = B(\U)\U_x.
\end{equation}
In the following, we first introduce a BP path-conservative method based on the flux globalization for \eqref{EQ:Hyp-NC-prod}. In order to apply the BP method proposed in this paper, we rewrite \eqref{EQ:Hyp-NC-prod} into a pseudo-conservative form
\begin{equation*}
\U_t+ \K(\U)_x = \mathbf{0},
\end{equation*}
where $\K$ is the global flux defined as follows
\begin{equation*}
\K(\U) = \F(\U) - \R(\U), \quad
\R(\U) = \int_{\widehat{x}}^x B(\U(\xi,t))\U_\xi(\xi,t)\,\mathrm{d}\xi.
\end{equation*}
The BP method, as studied in \S\ref{Sec:BP-FVM-Scalar}, can be applied once the values $\K_\jph^\pm$ and $\K_j$ are available.
To this end, we adopt the flux-globalization-based path-conservative method in \cite{KURGANOV2023111773} to compute $\K_\jph^\pm$ and $\K_j$ by
\begin{equation*}
\K_\jph^\pm = \F(\U_\jph^\pm) - \R_\jph^\pm, \quad
\K_j = \F(\U_j) - \R_j,
\end{equation*}
where $\U_j$ is a vector of values of the reconstruction polynomials of $\U$ at $x_j$. The terms $\R_\jph^\pm$ and $\R_j$ are obtained as follows. First, select $\widehat{x} = x_\half$, such that
\begin{equation*}
\R_\half^- = \mathbf{0}, \quad
\R_\half^+ = \B_{\bPsi,\half}.
\end{equation*}
Then, recursively compute $\R_\jph^\pm$ as
\begin{equation*}
\R_j = \R_\jmh^+ + \B_{j,L},\quad
\R_\jph^- = \R_j + \B_{j,R}, \quad
\R_\jph^+ = \R_\jph^- + \B_{\bPsi,\jph}.
\end{equation*}
The term $\B_j$ is the integration of the non-conservative product over left and right halves of $C_j$, that is,
\begin{equation*}
\B_{j,L} = \int_{x_\jmh}^{x_j} B(\U)\U_x \,\mathrm{d}x, \quad
\B_{j,R} = \int_{x_j}^{x_\jph} B(\U)\U_x \,\mathrm{d}x.
\end{equation*}
To evaluate the integrals, we first evaluate values at Newton-Cotes quadrature points $x_\jmh$, $x_j$, and $x_\jph$, using WENO polynomial introduced in \S\ref{Sec:WENO}. Then we recover a quadratic polynomial for each component involved in the computation based on these values. Finally, the non-conservative product is integrated on each half of the cell to obtain $\B_{j,L}$ and $\B_{j,R}$.
The terms $\B_{\bPsi,\jph}$ evaluates contributions of the non-conservative products across the cell interface $x_\jph$,
\begin{equation*}
\B_{\bPsi,\jph} = \int_0^1 B(\bPsi_\jph(s))\bPsi_\jph '(s)\,\mathrm{d}s,
\end{equation*}
where $\bPsi_\jph(s) := \bPsi(s; \U_\jph^-, \U_\jph^+)$ is a sufficiently smooth path connecting $\U_\jph^-$ and $\U_\jph^+$, that is,
\begin{equation*}
\bPsi: [0, 1] \times \mathbb{R}^d \times \mathbb{R}^d \to \mathbb{R}^d, \quad
\bPsi\left(0; \U_\jph^-, \U_\jph^+ \right) = \U_\jph^-, \quad
\bPsi\left(1; \U_\jph^-, \U_\jph^+ \right) = \U_\jph^+.
\end{equation*}

\begin{remark}
{\rm In the first-order version, the solution is represented by its cell averages $\overline{\U}_j$, so that $\U_j = \overline{\U}_j$, $\U_x = \bm{0}$, and
\begin{equation*}
\R_\jmh^+ = \R_j = \R_\jph^-.
\end{equation*}
}
\end{remark}

We now consider the five-equation transport model of two-medium flows \cite{ALLAIRE2002577}. In this model, the conservative variables, flux vector, and non-conservative products read as
\begin{equation*}
\U =
\begin{bmatrix}
z_1 \rho_1\\
z_2 \rho_2\\
\rho u\\
E\\
z_1
\end{bmatrix},\quad
\F(\U) =
\begin{bmatrix}
z_1 \rho_1 u\\
z_2 \rho_2 u\\
\rho u^2 + p\\
(E+p)u\\
0
\end{bmatrix},\quad
B(\U) =
\begin{bmatrix}
0&&&&\\
&0&&&\\
&&0&&\\
&&&0&\\
&&&&-u
\end{bmatrix}.
\end{equation*}
Here, $\zrho{k}$ are the partial densities, whose sum $\rho=\zrho{1}+\zrho{2}$ is the density in mixture, $z_k$ are the volume fractions of components $k\in\{1,2\}$, satisfying $z_1+z_2=1$, $u$ and $p$ are the velocity and pressure in mixture, respectively. The total energy $E$ consists of internal energy $\rho e$ and kinetic energy $\half\rho u^2$, that is,
\begin{equation*}
E=\rho e+\half\rho u^2.
\end{equation*}
The system is closed by the stiffened gas equation of state:
\begin{equation*}
p=(\gamma-1)\rho e-\gamma\pi_\infty,
\end{equation*}
with $\gamma$ and $\pi_\infty$ being the specific heat ratio and reference pressure of the mixture defined by
\begin{equation*}
\frac{1}{\gamma-1}=\frac{z_1}{\gamma_1-1}+\frac{z_2}{\gamma_2-1},\qquad
\frac{\gamma\pi_\infty}{\gamma-1}=\frac{z_1\gamma_1\pi_{\infty,1}}{\gamma_1-1}+\frac{z_2\gamma_2\pi_{\infty,2}}{\gamma_2-1}.
\end{equation*}
Here, the constants $\gamma_k$ and $\pi_{\infty,k}$, are the specific heat ratio and reference pressure of the component $k$. Moreover, the speed of sound reads as
\begin{equation*}
c = \sqrt{\frac{\gamma(p+\pi_\infty)}{\rho}}.
\end{equation*}
For the Five-equation model, the convex physics-constraint admissible set is
\begin{equation*}
\mG = \{\U \,\vert\, 0\leq z_1 \leq 1, \, \rho>0, \, \rho e -\pi_\infty>0 \},
\end{equation*}
which is convex if $(\gamma_1-\gamma_2)(\pi_{\infty,1}-\pi_{\infty,2})\geq 0$, which is satisfied by most gas-water interaction problems.

It is stressed that a proper selection of paths is crucial for the BP CFL conditions. For the five-equation model, paths for velocity $u$ and volume fraction $z_1$ need to be designed to evaluate $\B_{\Psi,\jph} = \int_0^1 u(s)z_1'(s)\,\mathrm{d}s$.
Inspired by \cite{Fu2025Bound}, we first introduce
\begin{equation*}
\alpha = \max_j\{ \dvert{u_j} + \widetilde{c}_j \},
\end{equation*}
where
\begin{equation*}
\widetilde{c}_j = \sqrt{c_j^2 + \kappa_j\frac{(\pi_\infty)_j}{\rho_j}}, \quad
\kappa_j = \begin{cases}
0, & p\geq 0,\\
1, & p < 0.
\end{cases}
\end{equation*}
We then define the path of velocity implicitly as
\begin{equation*}
\int_0^1 u_\jph(s)\,\mathrm{d}s = u_\jph^* = \frac{(\rho u)_\jph^*}{(z_1\rho_1)_\jph^* + (z_2\rho_2)_\jph^*}, \quad
\U_\jph^* = \frac{1}{2}\left(\U_\jph^- + \U_\jph^+\right) - \frac{1}{2\alpha}\left(\F(\U_\jph^+) - \F(\U_\jph^-)\right),
\end{equation*}
and use a linear path for $z_1$
\begin{equation*}
(z_1)_\jph(s) = (z_1)_\jph^- + s\left[(z_1)_\jph^+ - (z_1)_\jph^-\right],
\end{equation*}
such that $(z_1)'_\jph(s) = (z_1)_\jph^+ - (z_1)_\jph^-$ recovers the jump of $z_1$ across the cell interface.
The above selection of the path of the velocity is consistent with the discretization of the velocity in the mass and momentum equations when LF flux is applied. Besides, it can be proved that $\alpha \geq \dvert{u_\jph^*}$; consult \cite[Lemma 3.1]{Fu2025Bound} for the proof.

Based on the preparations above, we apply the modified LF flux on the moving mesh:
\begin{equation*}
\K_\jph -\omega_\jph \U_\jph \approx \mbH\left(\omega_\jph,\U_\jph^-,\U_\jph^+\right) = \half\left( \K_\jph^- + \K_\jph^+ \right) - \frac{\alpha}{2} \left( \U_\jph^+ - \U_\jph^- \right) - \omega_\jph \widehat{\U}_\jph^*,
\end{equation*}
where
\begin{equation}\label{EQ:Uast-5EQ}
\widehat{\U}_\jph^* = \frac{1}{2}\left(\U_\jph^- + \U_\jph^+\right) - \frac{1}{2\alpha}\left(\K(\U_\jph^+) - \K(\U_\jph^-)\right).
\end{equation}
The modification above ensures consistency of the velocity in the discretization, and thus leads to moderate BP CFL conditions.

Using the proposed BP method, we seek for the BP CFL conditions such that $\U_j^{(\ell),\mathrm{L},\pm} \in \mG$, where
\begin{equation}\label{EQ:split-5EQ}
\begin{aligned}
\U_j^{(\ell),\mathrm{L},-} &= \overline{\U}_j^{(\ell-1)} +  2\lambda_j^{(\ell-1)}\left[\mbH
\left(\omega_\jmh,\overline{\U}_{j-1}^{(\ell-1)},\overline{\U}_{j}^{(\ell-1)}\right) - \mbH\left(\omega_\jmh,\overline{\U}_j^{(\ell-1)},\overline{\U}_j^{(\ell-1)}\right)\right] \\
&= \left[ 1-2\lambda_j^{\ell-1}(\alpha^{(\ell-1)}-\omega_\jmh)\right]\overline{\U}_j^{(\ell-1)} + 2\lambda_j^{\ell-1}(\alpha^{(\ell-1)}-\omega_\jmh)\widehat{\U}_\jmh^{(\ell-1),*}, \\
\U_j^{(\ell),\mathrm{L},+} &= \overline{\U}_j^{(\ell-1)} - 2\lambda_j^{(\ell-1)}\left[\mbH
\left(\omega_\jph,\overline{\U}_{j}^{(\ell-1)},\overline{\U}_{j+1}^{(\ell-1)}\right) - \mbH\left(\omega_\jph,\overline{\U}_j^{(\ell-1)},\overline{\U}_j^{(\ell-1)}\right)\right] \\
&= \left[ 1-2\lambda_j^{\ell-1}(\alpha^{(\ell-1)}+\omega_\jph)\right]\overline{\U}_j^{(\ell-1)} + 2\lambda_j^{\ell-1}(\alpha^{(\ell-1)}+\omega_\jph)\widehat{\U}_\jph^{(\ell-1),*},
\end{aligned}
\end{equation}
where $\widehat{\U}_\jph^{(\ell-1),*}$ is defined in \eqref{EQ:Uast-5EQ}.
The BP conditions are stated and shown in Theorem \ref{Th:BP-CFL-5EQ} below.

\begin{lemma}\label{Lem:BP-CFL-5EQ}
Assume that $\overline{\U}_j^{(\ell-1)}\in\mG$ for all $j$. Then $\widehat{\U}_\jph^{(\ell-1),*}\in\mG$, where $\widehat{\U}_\jph^{(\ell-1),*}$ is given by (\ref{EQ:Uast-5EQ}).
\end{lemma}

\begin{proof}
Using the definition of $\K$, we rewrite $\widehat{\U}_\jph^{(\ell-1),*}$ as
\begin{equation*}
\begin{aligned}
\widehat{\U}_\jph^{(\ell-1),*} &= \left( (z_1\rho_1)_\jph^*, (z_2\rho_2)_\jph^*, (\rho u)_\jph^*, E_\jph^*, (\widehat{z}_1)_\jph^* \right)^\top, \\
(\widehat{z}_1)_\jph^* &= \frac{u_\jph^* + \alpha^{(\ell-1)}}{2\alpha^{(\ell-1)}} \overline{z_1}_j^{(\ell-1)} + \frac{\alpha^{(\ell-1)} - u_\jph^*}{2\alpha^{(\ell-1)}} \overline{z_1}_{j+1}^{(\ell-1)},
\end{aligned}
\end{equation*}
which is a simplified version of the central-upwind scheme in \cite{Fu2025Bound}, and thus $\widehat{\U}_\jph^{(\ell-1),*}\in\mG$; see detailed proof in \cite[Theorem 3.1]{Fu2025Bound}.
\end{proof}

\begin{theorem}\label{Th:BP-CFL-5EQ}
Assume that $\overline{\U}_j^{(\ell-1)}\in\mG$ for all $j$. Then $\U_j^{(\ell),\rm L,\pm}\in\mG$, where $\U_j^{(\ell),\rm L,\pm}$ is given by (\ref{EQ:split-5EQ}), provided that
\begin{enumerate}
\item the grid velocity $\omega_\jph$ satisfies $\alpha\geq \dvert{\omega_\jph}$, and
\item the BP CFL condition
\end{enumerate}
\begin{equation*}
\lambda_j^{(\ell-1)}\left(\alpha^{(\ell-1)}+\dvert{\omega_\jph}\right) = \frac{\dt^n \left(\alpha^{(\ell-1)}+\dvert{\omega_\jph}\right)}{\dx_j^{(\ell-1)}+\dt^n(\omega_\jph - \omega_\jmh)}
 \leq \frac{1}{2}
\end{equation*}
is satisfied.
\end{theorem}

\begin{proof}
We only prove that $\U_j^{(\ell),\rm L,+}\in\mG$, as $\U_j^{(\ell),\rm L,-}\in\mG$ can be proved in a similar way.
The splitting in \eqref{EQ:split-5EQ} is convex if the two conditions in Theorem \ref{Th:BP-CFL-5EQ} are satisfied. By Lemma \ref{Lem:BP-CFL-5EQ}, we have $\widehat{\U}_\jph^{(\ell-1),*}\in\mG$. It follows from the convexity of $\mG$ that $\U_j^{(\ell),\rm L,+}\in\mG$. This completes the proof the theorem \ref{Th:BP-CFL-5EQ}.
\end{proof}

According to Theorem \ref{Th:BP-CFL-5EQ}, the grid velocity should satisfy $\alpha\geq \dvert{\omega_\jph}$. To fulfill the requirement, we modify $\omega_\jph$ as
\begin{equation*}
\omega_\jph^{(\mathrm{mod})} = \begin{cases}
\omega_\jph, & \dvert{\omega_\jph} \leq \alpha, \vspace{+3pt}\\
\operatorname{sign}(\omega_\jph)\,\alpha,  &\dvert{\omega_\jph} > \alpha,
\end{cases}\qquad
\operatorname{sign}(\phi) = \begin{cases}
 1, & \phi\geq 0,\\
-1, & \phi < 0.
\end{cases}
\end{equation*}

We point out that the BP CFL condition in Theorem \ref{Th:BP-CFL-5EQ}  is not severe, as the velocity $u_\jph^*$ used in the fifth equation is consistent with the velocity hidden in the discretizations of mass and momentum equations when the LF numerical flux is applied; consult also \cite{Fu2025Bound}.

We also discuss a few implementation issues here. Firstly,
similar to \S\ref{Sec:BP-CFL}, one can modify high-order sub-cell schemes $\U_j^{(\ell),\mathrm{H},\pm}$ by blending them with their first-order counterparts $\U_j^{(\ell),\mathrm{L},\pm}$ the same way as \eqref{EQ:blend-RKl}, or modify high order numerical flux $\K_\jph^{(\ell-1)}$ the same way as \eqref{EQ:flux-RKl-blending}. Moreover, similar to the Euler equations, the BP flux limiters are first employed to enforce positivity of partial densities ($z_1\rho_1$ and $z_2\rho_2$) and bound of volume fraction ($z_1$), and are then applied to preserve positivity of $\rho e -\pi_\infty$.

\section{Numerical examples}\label{Sec:Num}

In this section, we will demonstrate robustness and effectiveness of our BP adaptive moving mesh method by a number of numerical examples. For simplicity, we denote methods with and without BP limiters by $\textbf{BP}$ and $\textbf{NBP}$, respectively. In addition, the moving mesh and uniform mesh are denoted by $\textbf{MM}$ and $\textbf{UM}$, respectively. In all numerical experiments, we set CFL = 0.16, which satisfies the AP CFL condition \eqref{EQ:AP-CFL}, as well as the BP CFL condition \eqref{EQ:BP-CFL}.

For moving mesh methods, we use a commonly used monitor function $\sigma$ (see Eq. (\ref{EQ:Euler-Lagrange})) based on \eqref{EQ:sigma}.
The convergence rate is computed by the following standard approach:
\begin{equation*}
\mathcal{R}^{2N} = \dfrac{\log_{10} \left( \varepsilon^{N}/\varepsilon^{2N} \right)}{\log_{10} \left( \dx_{\max}^{N}/\dx_{\max}^{2N} \right)}.
\end{equation*}

\begin{figure}[ht!]
\centering	
\includegraphics[width=0.48\textwidth]{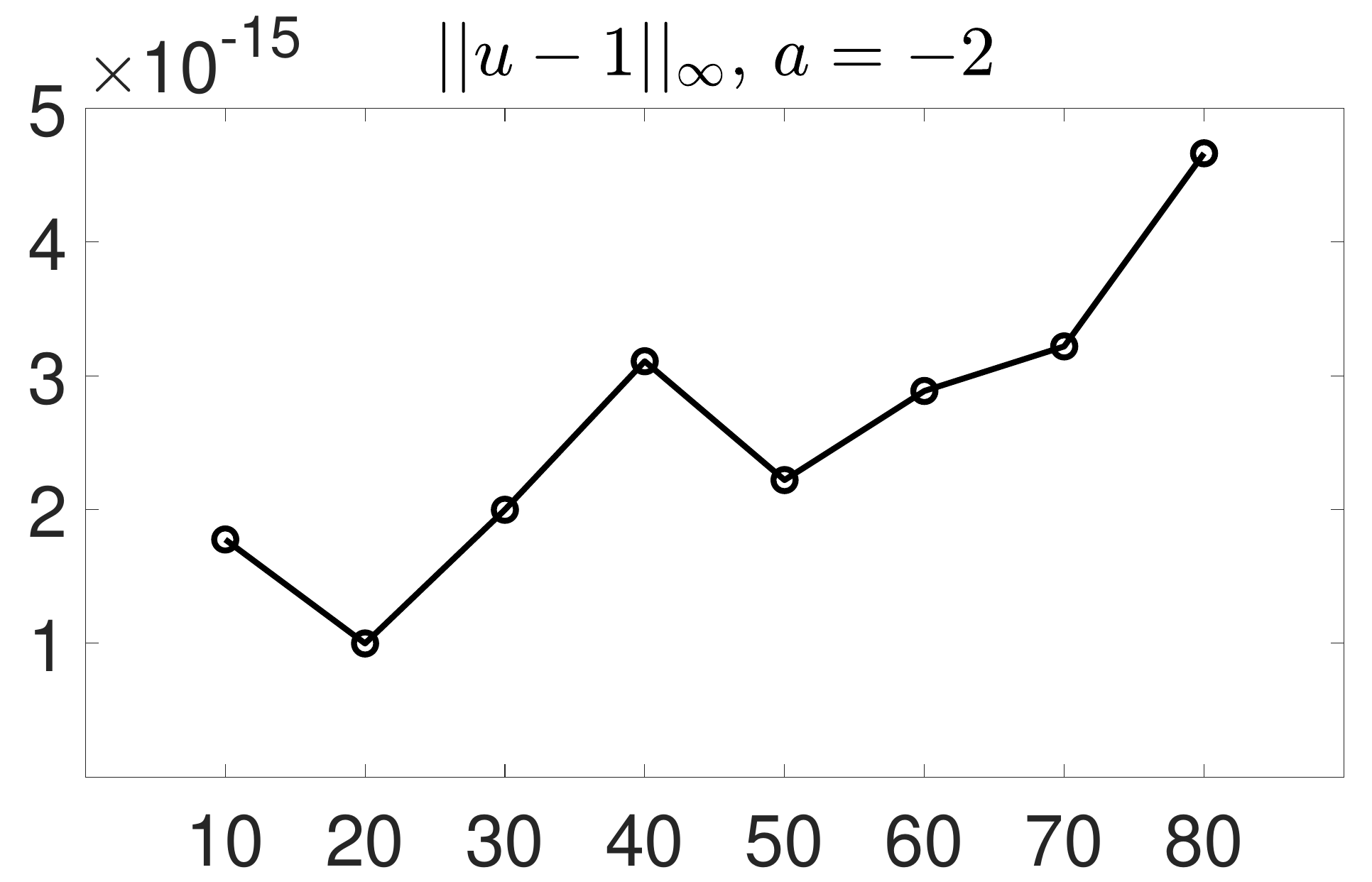}\qquad
\includegraphics[width=0.48\textwidth]{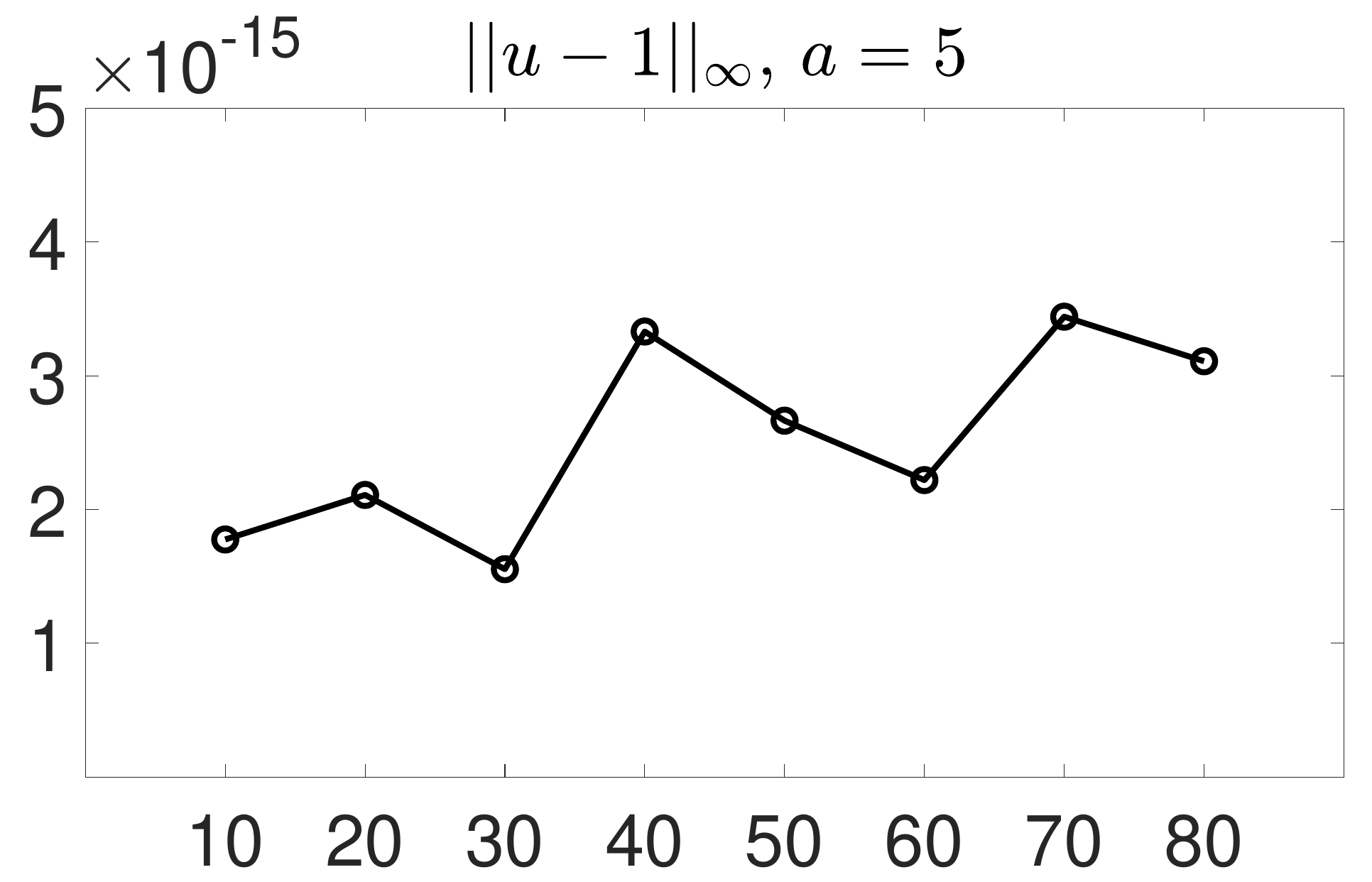}
\caption{\ttfamily Example \ref{exam-GCL}. D-GCL test. Left: $a=-2$. Right: $a=5$. The output time is $t=2$.\label{Fig:GCL}}
\end{figure}

\begin{example}\label{exam-GCL} 
{\em advection equation.}
In the first example, the D-GCL is verified by simulating the linear advection equation
\begin{equation*}
u_t + a u_x = 0, \quad
u(x, 0) = 1, \quad
x\in [0,2\pi],
\end{equation*}
which admits a uniform solution $u(x,t) = 1$. Periodic boundary conditions are employed on the boundaries of the physical domain.
\end{example}
In this example, we initialize $Du$ in each control volume by an arbitrary random number between $0$ and $1$. Besides, we set $\beta=0.6$. The $L^\infty$ errors $\dVert{u-1}_\infty$ at $t=2$ are shown in Figure \ref{Fig:GCL} with different $a$ and $N$. 
One can see that the errors are close to machine zero, which indicates that the D-GCL is satisfied by our approach.

\begin{example}\label{exam1}
{\em Burgers equation.} 
In the second test, we verify BP and AP properties of the proposed method by solving the nonlinear inviscid Burgers equation
\begin{equation*}
u_t + (\frac{1}{2}u^2)_x = 0, \quad 0< x < 2\pi,
\end{equation*}
subject to periodic boundary conditions and the initial condition $u_0(x)=\sin^4(x)$. 
\end{example}

The invariant domain for this test is $\mG =\{u\vert\; 0\leq u\leq 1\}$.
In this example, we use the monitor function \eqref{EQ:sigma} with 
\begin{equation*}
Du = \sqrt{\dvert{u_x}^2 + \dvert{u_{xx}}^2}
\end{equation*}
with $\beta=0.6$. We first solve the problem until $t=0.4$ when the solution is still smooth. The exact solution is obtained using Newton iteration.  The $L^1$-error, convergence rate, and the minimum numerical solution are shown in Table \ref{Tab:Accuracy-Bu-CFL16}. It can be observed that the BP method preserves the bound of the solution on moving meshes, without sacrificing the accuracy. Therefore, the AP and BP properties of the scheme are verified by the test. However,  although the conventional moving mesh methods without BP limiters achieve third-order of convergence rate, the numerical solution exceeds the bound (between $0$ and $1$).

\begin{table}[htb!]
\caption{\ttfamily Convergence of AMM-FVM schemes with/without BP limiters, together with minimum of solution, for Burgers equation. The $\CFL$ number is taken as $0.16$.}
\label{Tab:Accuracy-Bu-CFL16}
\begin{center}
\begin{tabular}{||c|cccc|cccc||}
\hline
\multirow{2}{*}{$N$}
&\multicolumn{4}{c}{Without BP Limiters}
&\multicolumn{4}{c||}{With BP Limiters} \\ \cline{2-9}
&$L^1$-error &$\dx_{\max}$ &Rate $\mathcal{R}$ &$u_{\min}$
&$L^1$-error &$\dx_{\max}$ &Rate $\mathcal{R}$ &$u_{\min}$\\ \hline
$ 40$ &1.27E-02 &2.33E-01 &--   &-1.11E-03
      &6.85E-03 &2.02E-01 &--   & 1.98E-08 \\

$ 80$ &3.68E-03 &1.53E-01 &2.95 &-1.90E-04
      &1.63E-03 &1.31E-01 &3.32 & 2.32E-10 \\

$160$ &7.45E-04 &9.51E-02 &3.36 &-6.00E-05
      &3.60E-04 &8.27E-02 &3.28 & 1.00E-16 \\
\hline
\end{tabular}
\end{center}
\end{table}

\begin{figure}[ht!]
\centering	
\includegraphics[width=0.32\textwidth]{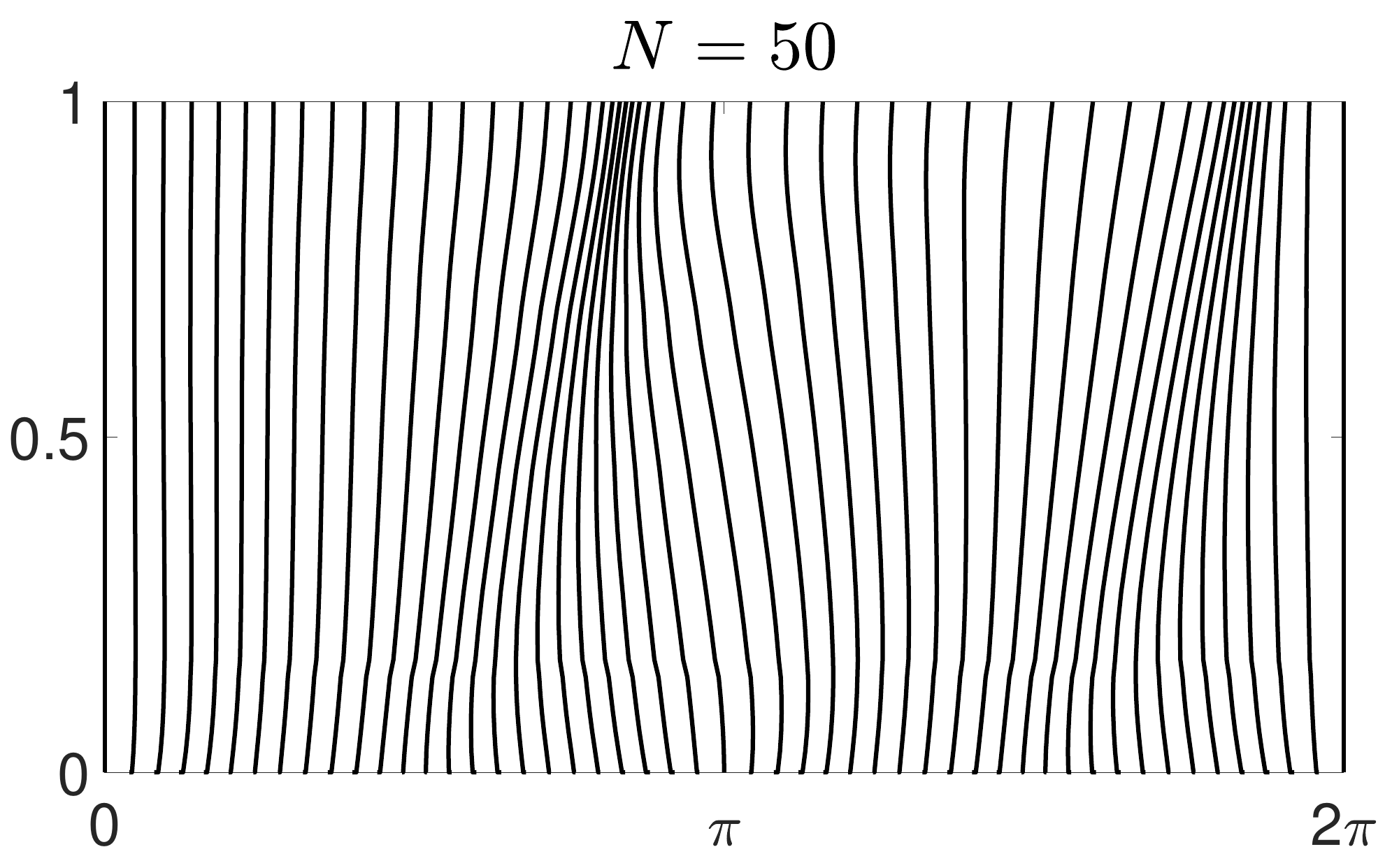}
\hspace*{0.01cm}
\includegraphics[width=0.32\textwidth]{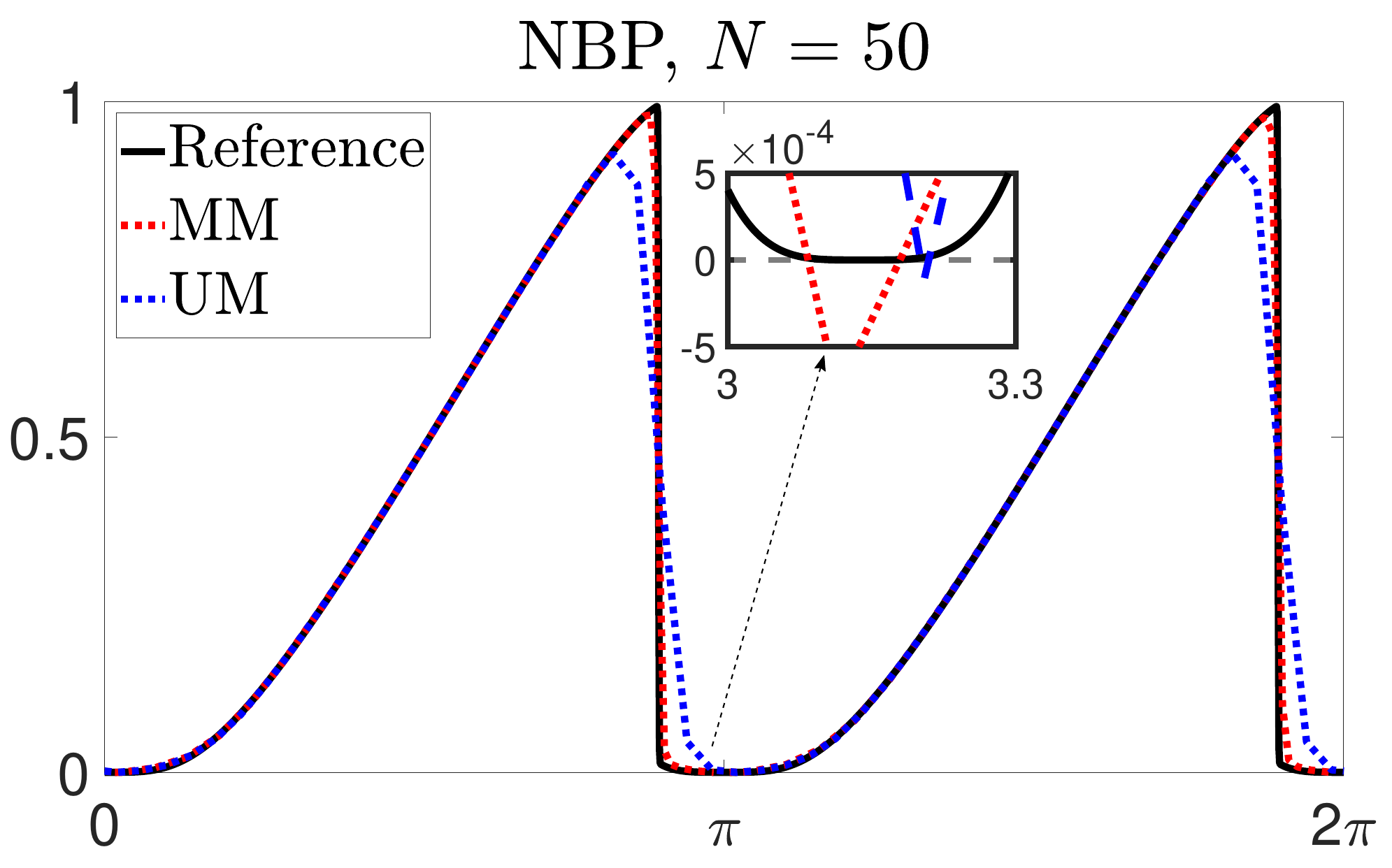}
\hspace*{0.01cm}
\includegraphics[width=0.32\textwidth]{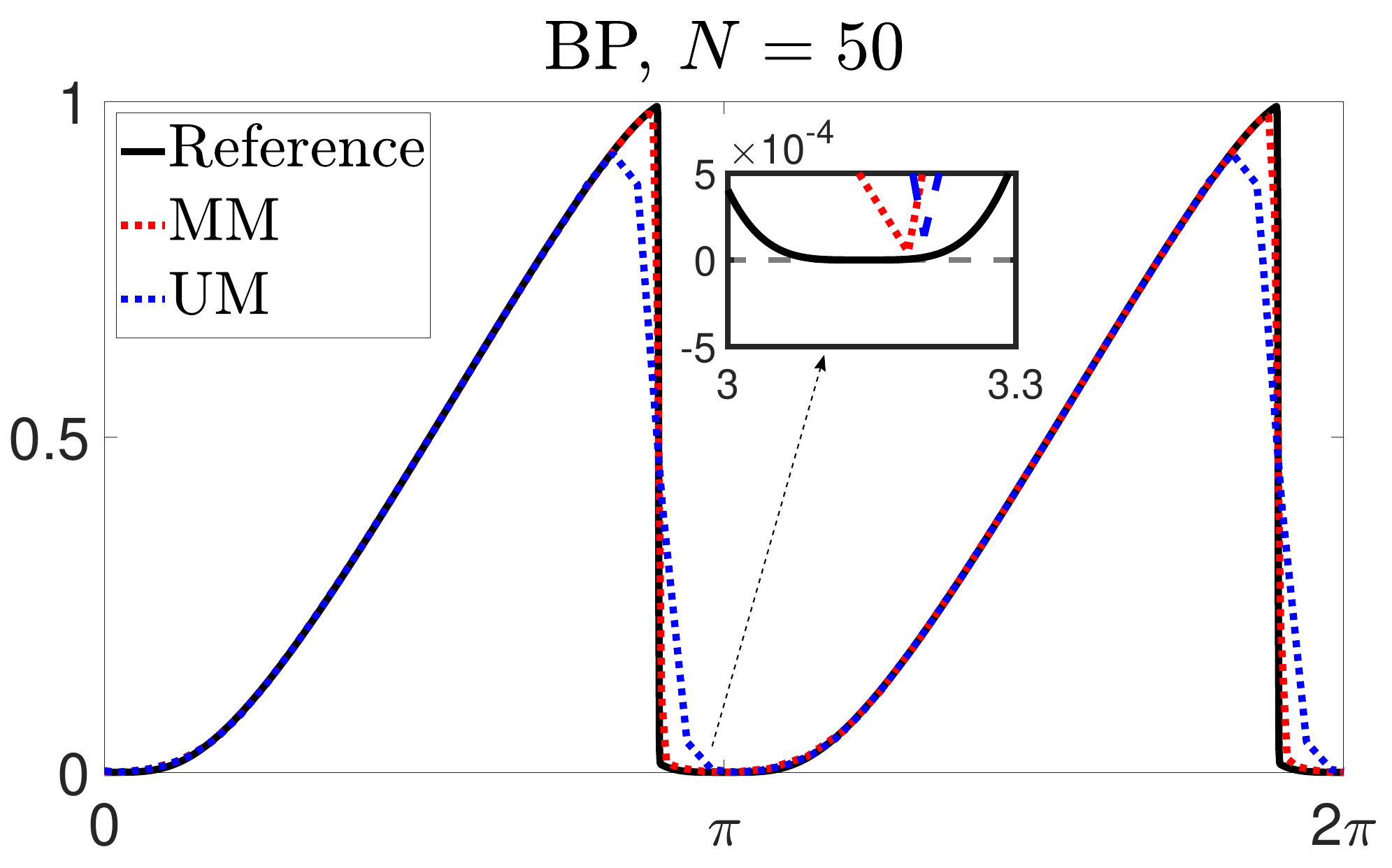}
\caption{\ttfamily Example \ref{exam1}.  Left: trajectory of grid points. Middle: comparison of solutions on uniform and adaptive meshes without BP limiters. Right: comparison of solutions on uniform and adaptive meshes with BP limiters. The output time is $t=1.3$.\label{Fig:Burgers}}
\end{figure}

Next, we simulate the problem with $N=50$ until $t=1.3$, when two shocks form. We use the same monitor function with $\beta=0.6$. The trajectory of grid points and comparison of numerical solutions are presented in Figure \ref{Fig:Burgers}. We compare the results with and without adaptive mesh redistribution. By using the mesh adaptivity, the mesh grids adaptively move towards shock waves, as illustrated in the first panel of Figure \ref{Fig:Burgers}. As a result, the numerical solutions on adaptive mesh (both NBP and BP) have superior resolution than those without mesh adaptivity; see the second and the third panels. In addition, the bound of numerical solutions is preserved both on uniform and adaptive meshes by using the proposed BP method (the third panel), whereas it is not the case for NBP finite-volume schemes (the second panel).

\begin{example}\label{exam2} 
{\em Euler equations.} 
In this test, we consider the Euler equations with the initial condition:
\begin{equation*}
(\rho, u, p)\Big|_{(x,0)} = \begin{cases}
(2, 0, 10^6), &x<0,\\
(1, 0, 1), & x\geq 0.
\end{cases}
\end{equation*}
Free boundary conditions are employed at the boundaries of the physical domain $[-1,1]$. The output time is $t=8\times 10^{-4}$.
\end{example}

The invariant domain for this test is $\mG =\{\U\mid \rho(\U)>0, \, p(\U)>0\}$.
In this example, we use the monitor function \eqref{EQ:sigma} with
\begin{equation*}
D\U = \sqrt{\dvert{\rho_x}^2 + \dvert{\rho_{xx}}^2},
\end{equation*}
with $\beta=0.3$. 

\begin{figure}[ht!]
\centering
\includegraphics[width=0.48\textwidth]{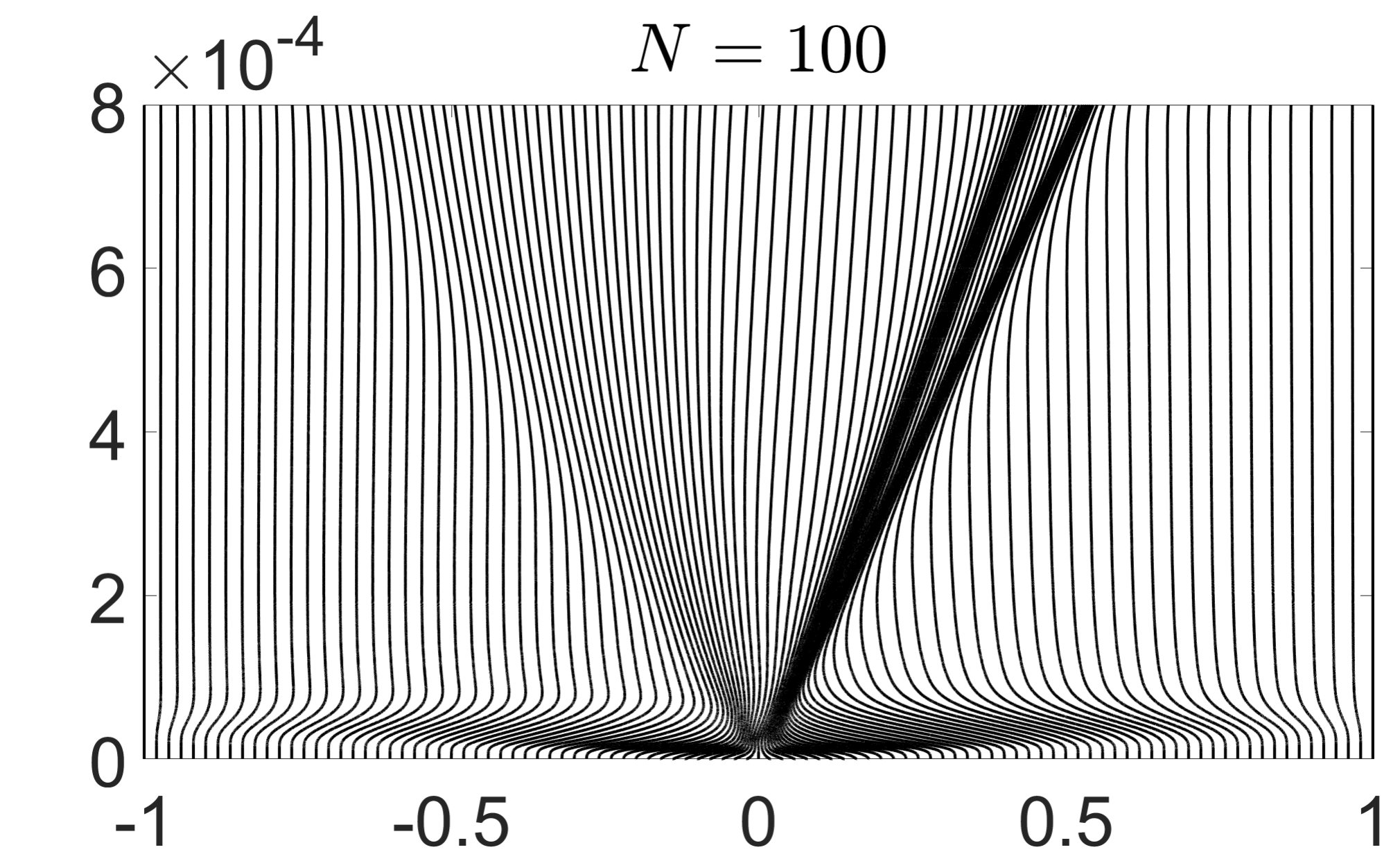}
\hspace*{0.05cm}
\includegraphics[width=0.48\textwidth]{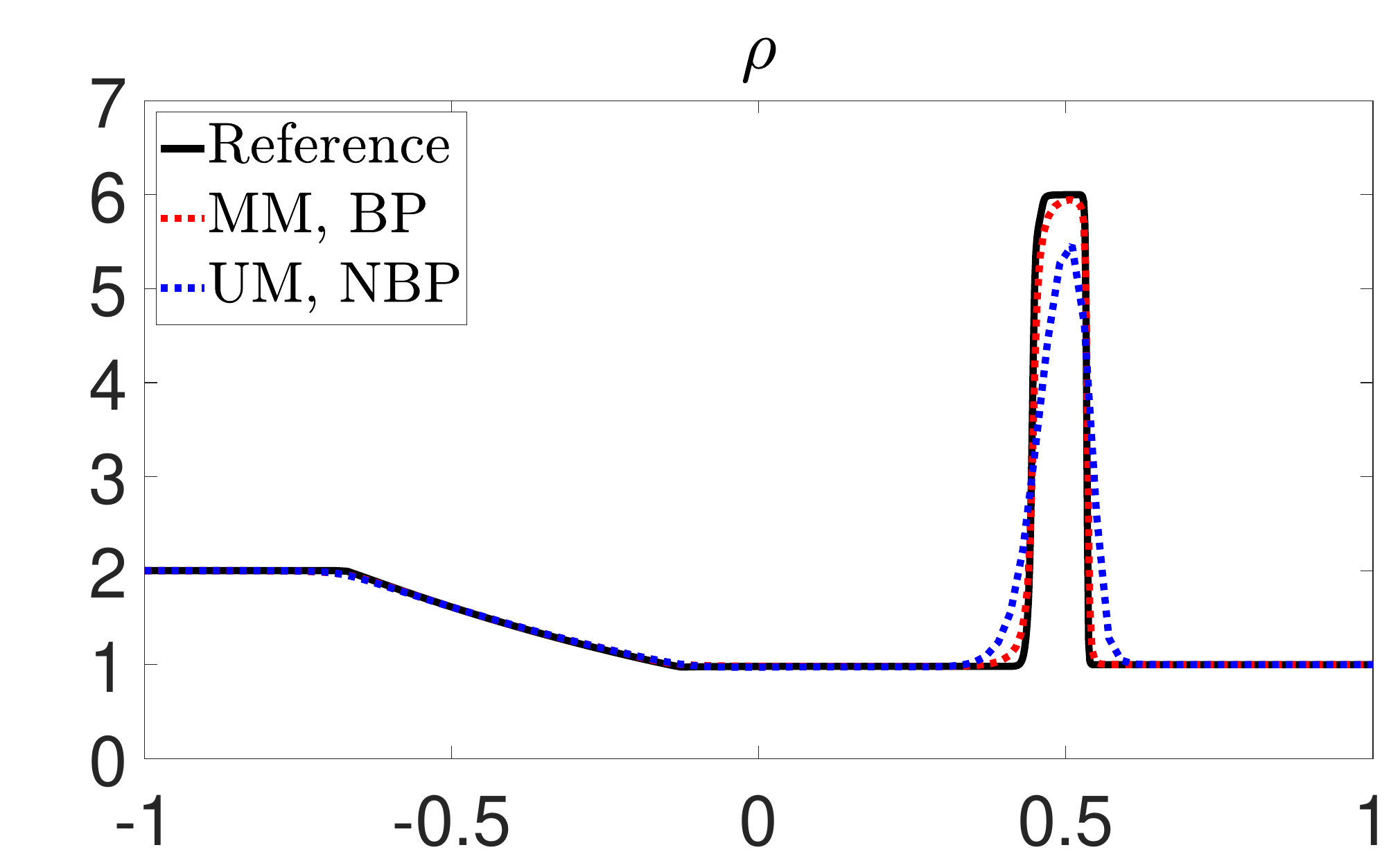}\\
\includegraphics[width=0.48\textwidth]{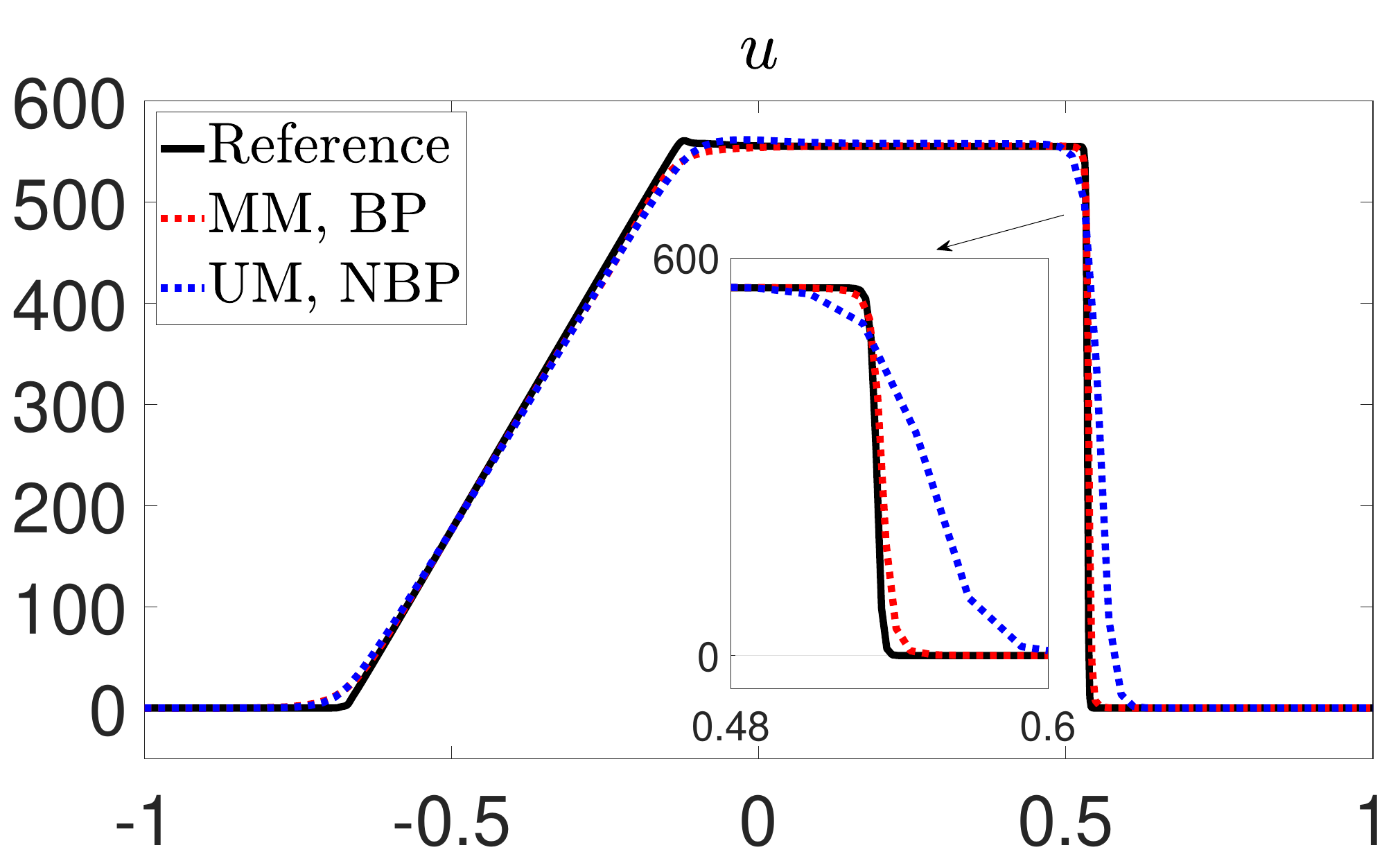}
\hspace*{0.05cm}
\includegraphics[width=0.48\textwidth]{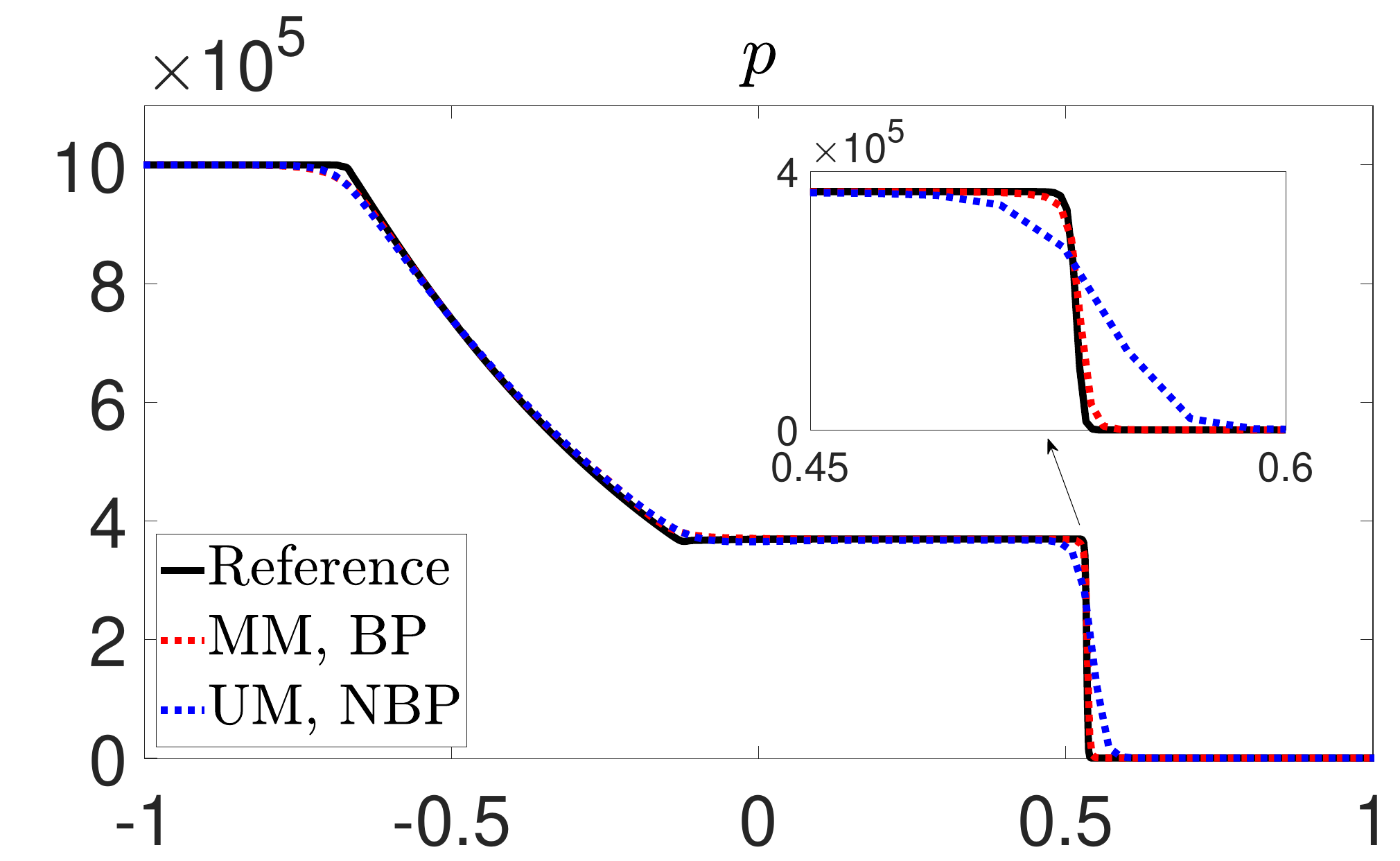}
\caption{\ttfamily Example \ref{exam2}: Trajectory of mesh grids (top left) and distributions of density $\rho$ (top right), velocity $u$ (bottom left), and pressure $p$ (bottom right) of the Riemann problem of Euler system at $t=8\times 10^{-4}$. \label{Fig:RP}}
\end{figure}

In this test, the code crashes quickly on adaptive mesh without utilizing BP limiters, due to large ratio in pressure. With the BP limiters proposed in this work, adaptive moving mesh results have demonstrated big advantages over the uniform mesh ones, which can be obtained without BP limiters. To see this, we plot in Figure \ref{Fig:RP} the trajectory of grid points and distributions of density, velocity, and pressure on uniform and adaptive meshes.  By using adaptive mesh redistribution, the rarefaction, contact and shock waves can be distinguished by the grid distribution against time evolution, as shown in the top left of Figure \ref{Fig:RP}. As a result, the density, velocity, and pressure are well approximated by using just $100$ cells, which is remarkably superior to those on uniform mesh with the same number of cells.

\begin{example}\label{exam3} 
{\em Five-equation model.} 
We consider the so-called five-equation model for two-medium flows. 
The gas-liquid shock-tube problem taken from \cite{CORALIC201495} will be considered. In this example, highly compressed air and water at atmospheric pressure form the left and right states, respectively. The initial data are
\begin{equation*}
(\rho_1, \rho_2, u, p, z_1)\Big|_{(x,0)} =
\begin{cases}
(1.241, 0.991, 0, 2.753, 1-10^{-13}), &-5\leq x<0,\\
(1.241, 0.991, 0, 3.059\times 10^{-4}, 10^{-13}), &0\leq x\leq 5,
\end{cases}
\end{equation*}
subject to free boundary conditions. Here, $\gamma = 1.4$, $\pi_{\infty,1} = 0$ and $\gamma_2 = 5.5$, $\pi_{\infty,2} = 1.505$. The output time is $t=1$.
\end{example}

In the simulation, the BP limiter is necessary for the volume fraction to be bounded on moving mesh. In this example, we use the monitor function with
\begin{equation} \label{mx1}
\text{MM(v1):} \quad
D\U = \sqrt{100\dvert{\rho_x}^2 + \dvert{\rho_{xx}}^2 + 1000\dvert{z_x}^2},
\end{equation}
and $\beta=0.3$. 

\begin{figure}[ht!]
\centering
\includegraphics[width=0.48\textwidth]{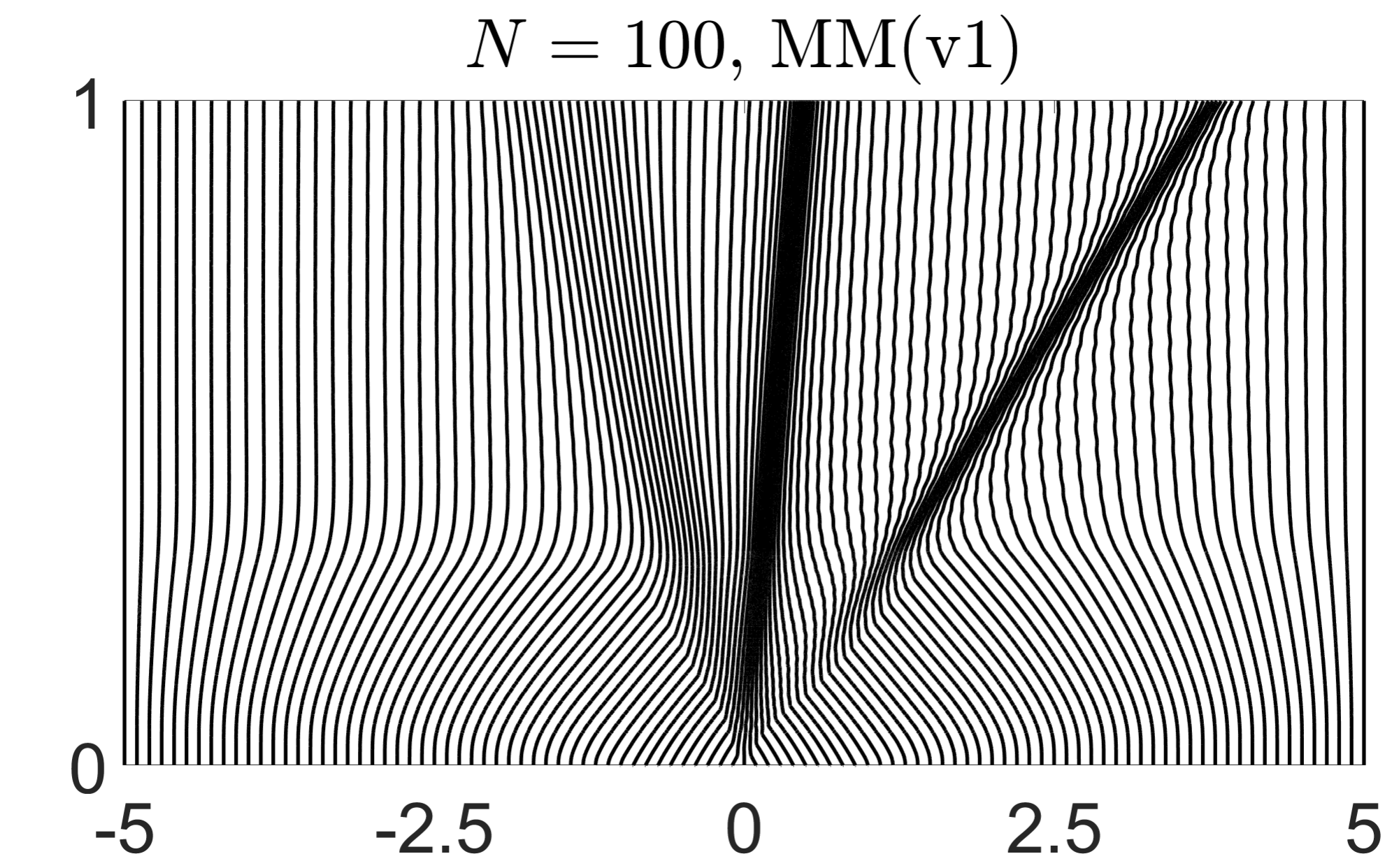}
\hspace*{0.05cm}
\includegraphics[width=0.48\textwidth]{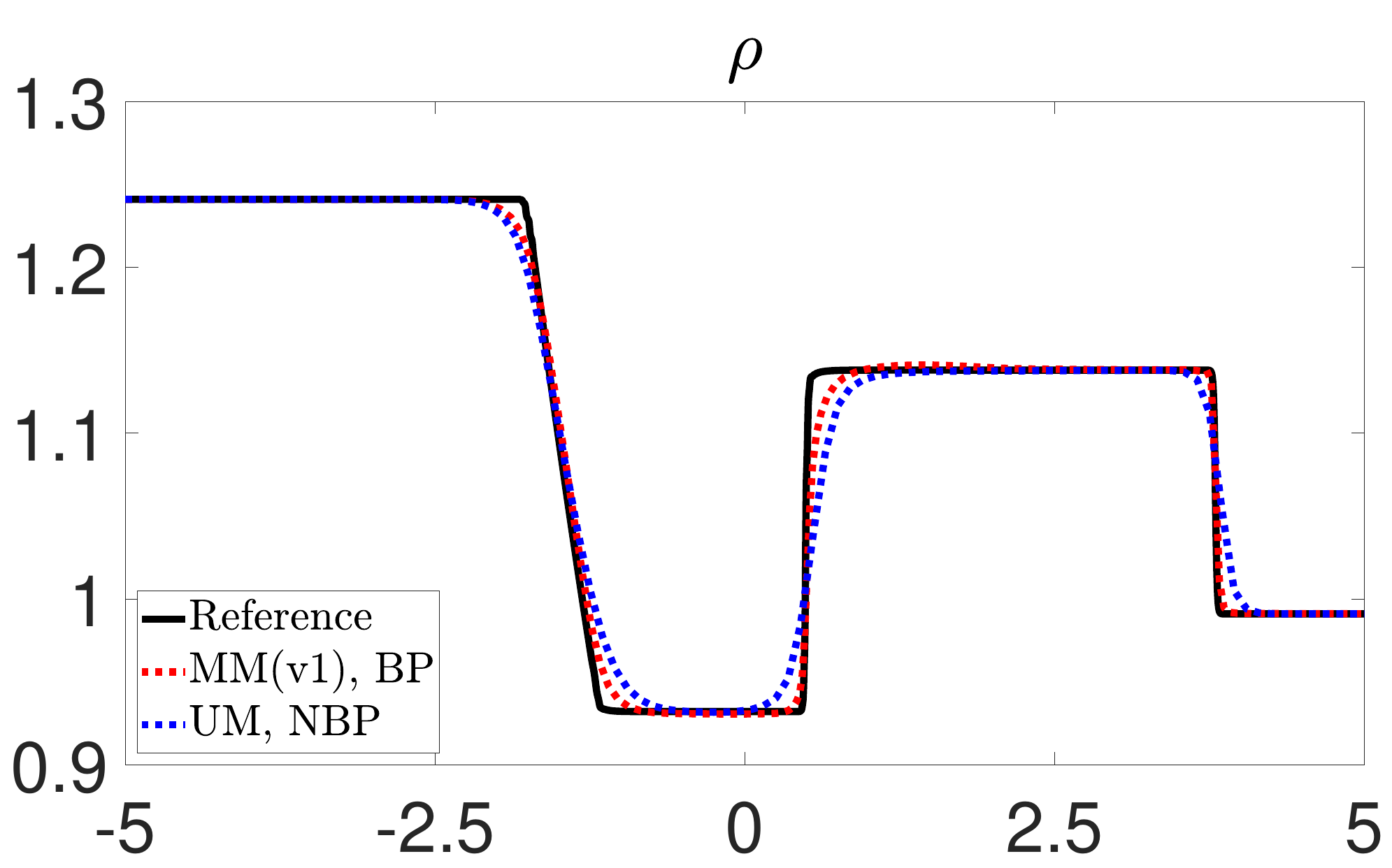}\\
\caption{\ttfamily Example \ref{exam3}: Trajectory of moving mesh grids (left) and distributions of density at $t=1$(right), where red color indicates the BP moving mesh solutions, and blue color presents the uniform mesh solutions. v1 represents results with the monitor function \eqref{mx1}. \label{Fig:5EQ-RP1}}
\end{figure}

The trajectory of grid points and a comparison of the density produced on adaptive and uniform meshes are illustrated in the first row of Figure \ref{Fig:5EQ-RP1}, where $100$ control volumes are used. It is found that the solution produced by moving mesh methods is much better than that produced on the uniform mesh, which indicates effectiveness of the proposed BP moving mesh methods. 

\begin{figure}[ht!]
\centering
\includegraphics[width=0.48\textwidth]{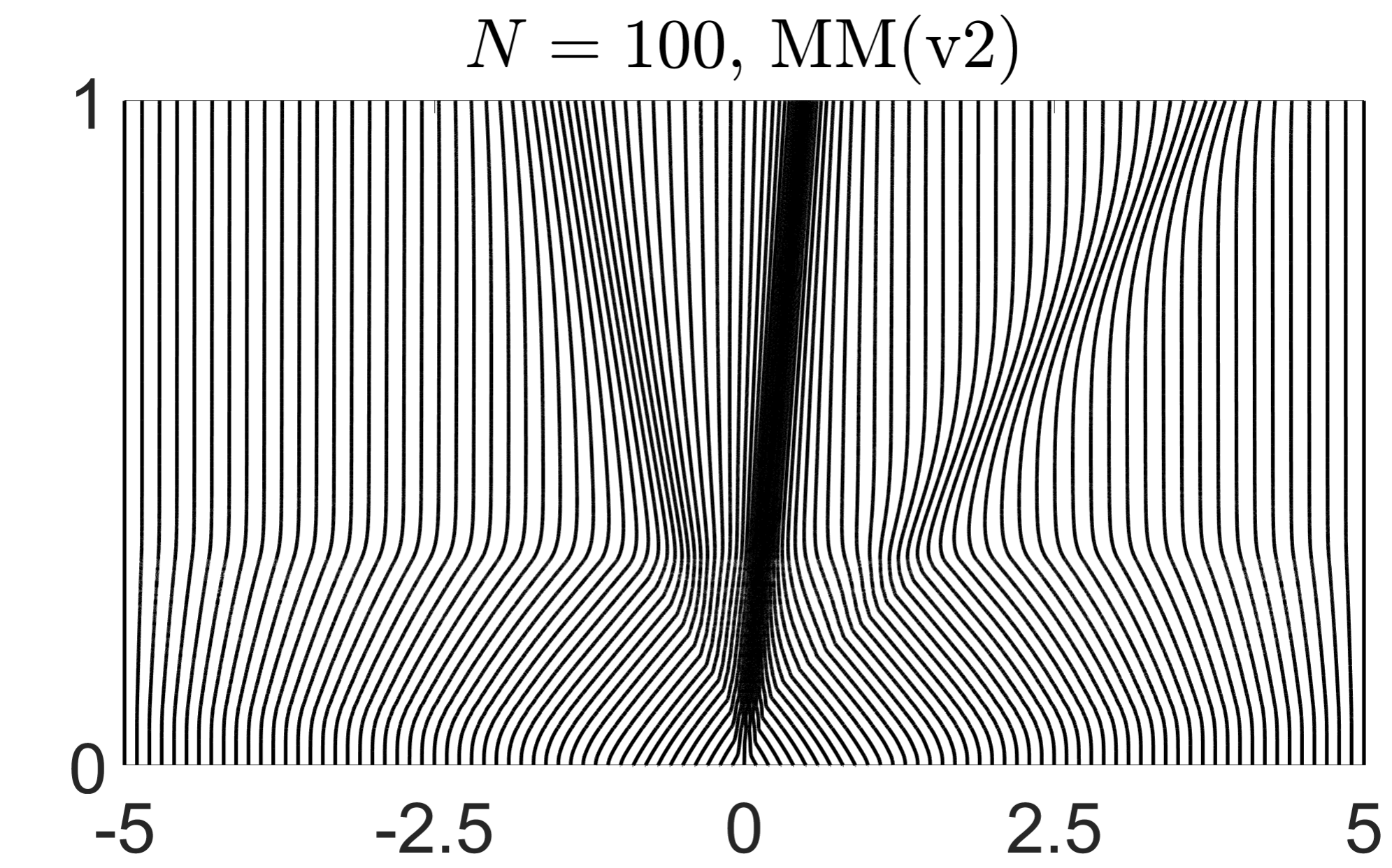}
\hspace*{0.05cm}
\includegraphics[width=0.48\textwidth]{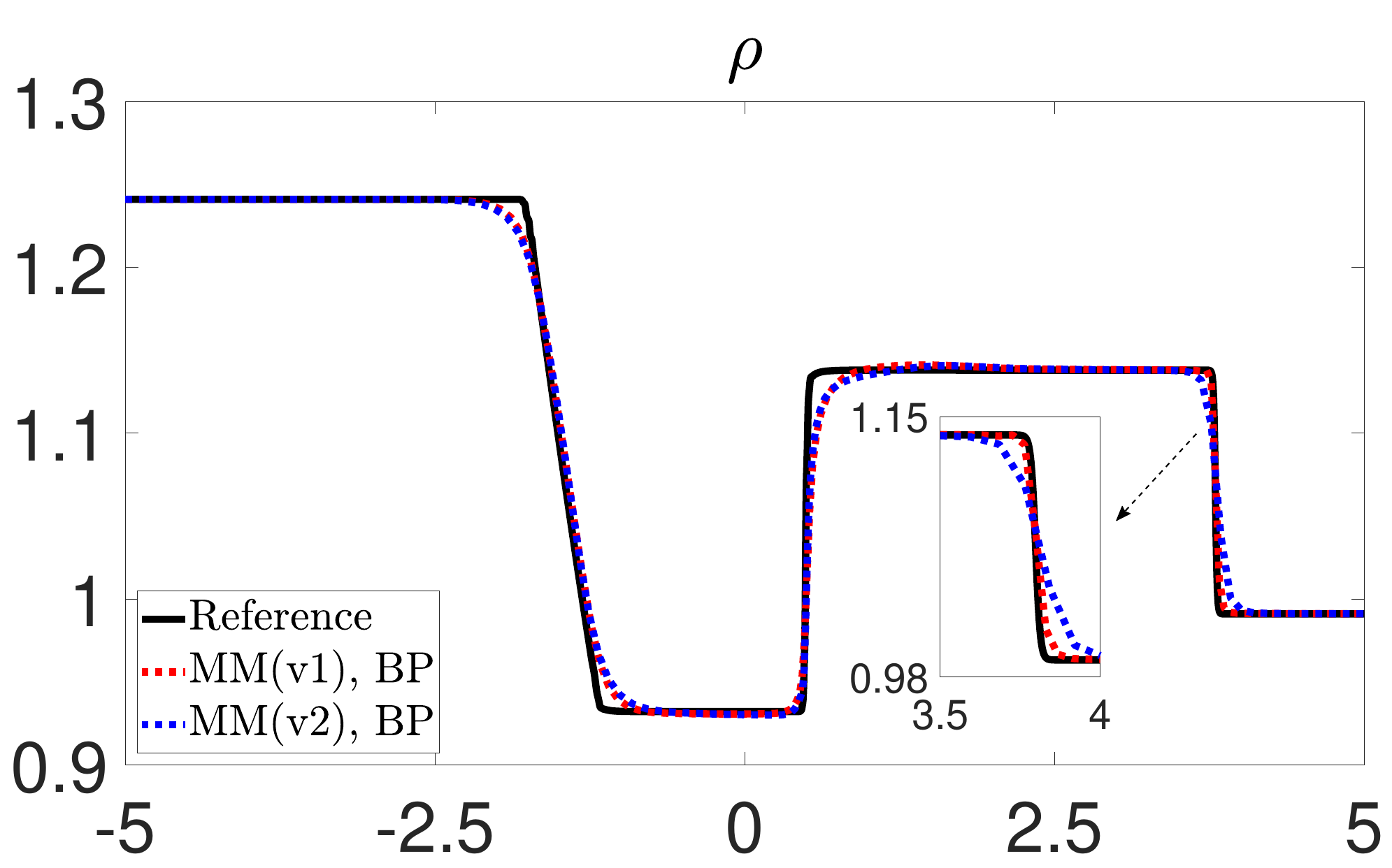}
\caption{\ttfamily Example \ref{exam3}: Trajectory of mesh grids (left) and distributions of density at $t=1$(right), v1 and v2 represent results with the monitor function \eqref{mx1} and  \eqref{mx2}, respectively.\label{Fig:5EQ-RP1A}}
\end{figure}

We point out that the second derivative term $\dvert{\rho_{xx}}^2$ in the monitor function (\ref{mx1}) plays an important role in improving quality of the mesh, which can detect the change of curvature effectively. As a comparison test, we only remove the term $\dvert{\rho_{xx}}^2$ in (\ref{mx1}), and use the monitor function with
\begin{equation} \label{mx2}
\text{MM(v2):} \quad
D\U = \sqrt{100\dvert{\rho_x}^2 + 1000\dvert{z_x}^2},
\end{equation}
and $\beta=0.3$. It can be seen from the left of Figure \ref{Fig:5EQ-RP1A} that the mesh trajectory can not well capture the rarefaction and shock waves, as compared with that in Figure \ref{Fig:5EQ-RP1}. As a result, the quality of the BP moving mesh result with monitor function v1 is remarkably
superior to that with v2, as can be seen in the right panel of Figure \ref{Fig:5EQ-RP1A}.

\begin{example} \label{exam4}
{\em Five-equation model.} 
This example tests a very challenging gas-water Riemann problem taken from \cite{CHENG2014143,Fu2025Bound}. Initially, the fluid to the left of the interface located at $x = 0.3$ is the ideal gas, while the fluid on the right is water. The initial conditions are
\begin{equation*}
(\rho_1, \rho_2, u, p, z_1)\Big|_{(x,0)} =
\begin{cases}
(5, 10^3, 0, 10^5, 1-10^{-13}), &x<0.3,\\
(5, 10^3, 0, 10^9, 10^{-13}), &x\geq 0.3,
\end{cases}
\end{equation*}
where $\gamma_1 = 1.4$, $\pi_{\infty,1} = 0$, and $\gamma_2 = 4.4$, $\pi_{\infty,2} = 6\times 10^8$. The problem is computed in $[0,1]$ with free boundary conditions. The final time is $2.4\times 10^{-4}$.
\end{example}

\begin{figure}[ht!]
\centering
\includegraphics[width=0.48\textwidth]{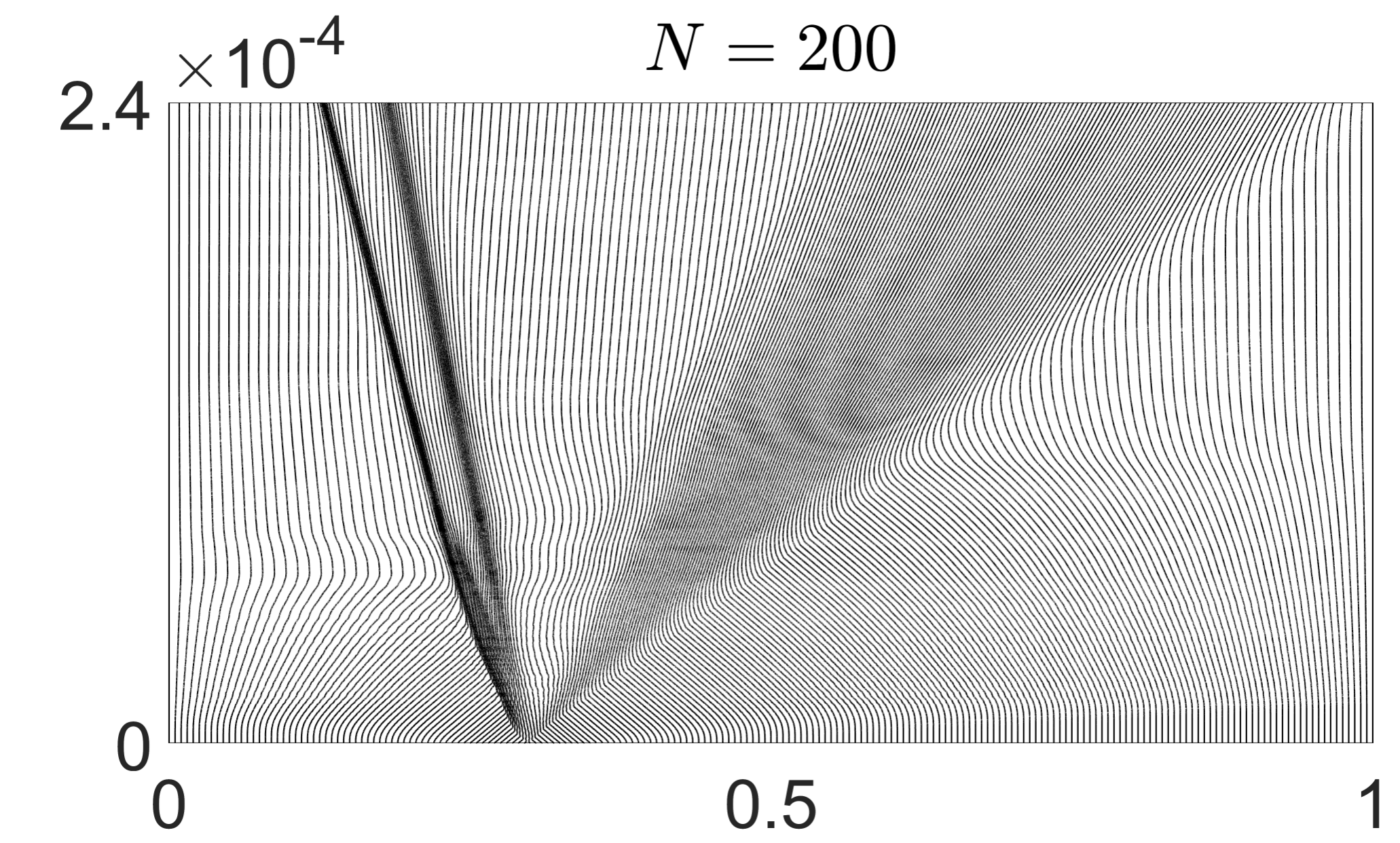}
\includegraphics[width=0.48\textwidth]{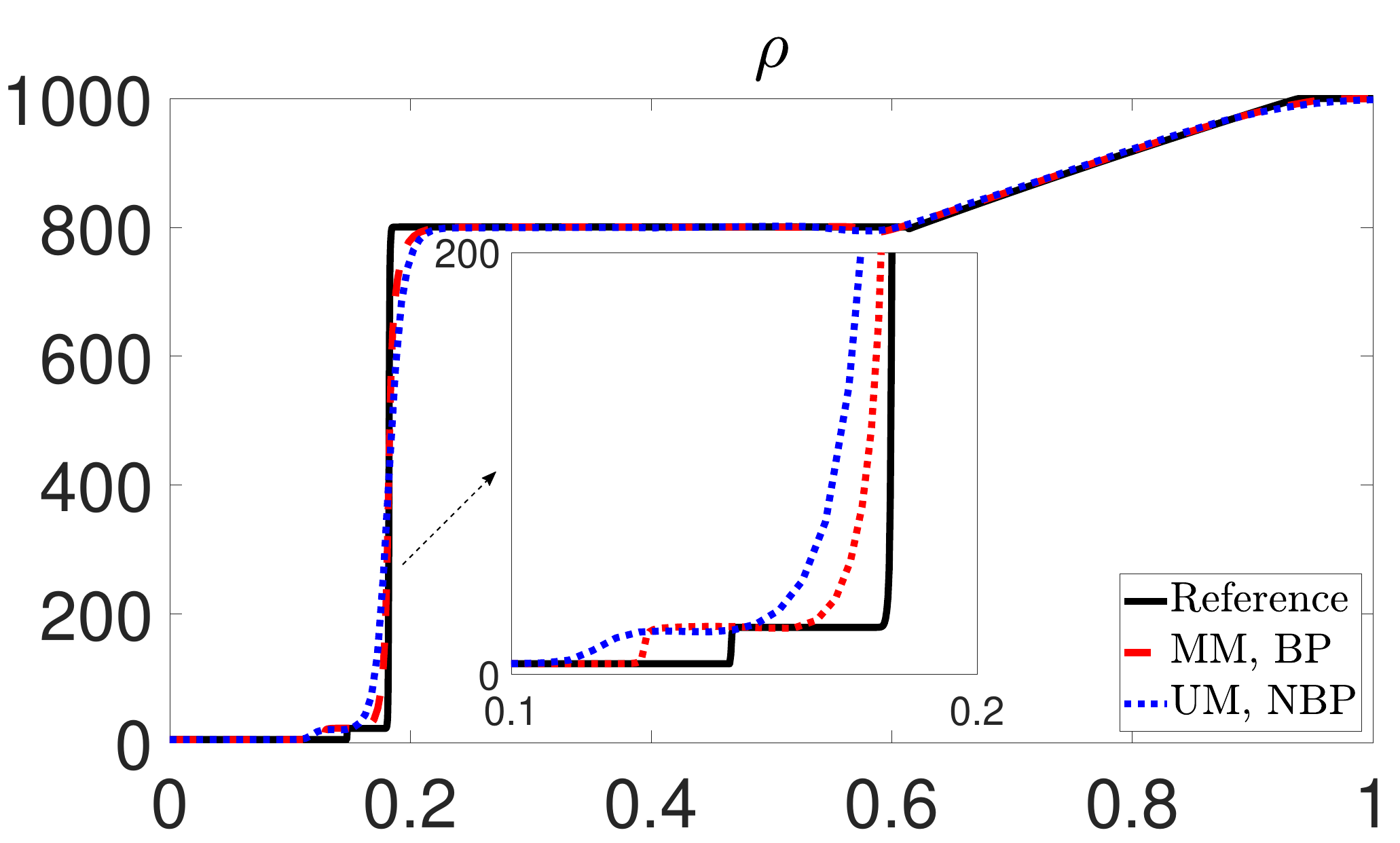}
\caption{\ttfamily Example \ref{exam4}: The trajectory of the moving mesh grids (left), and the distributions of density at $t=2.4\times 10^{-4}$ (right). \label{Fig:5EQ-RP2}}
\end{figure}

In this example, the amplitude of shock wave is much smaller than that of the material interface. Besides, the shock wave and the material interface are quite close to each other. These two factors make it difficult to capture both structures on uniform mesh. In addition, the large ratios in density and pressure make this problem challenging on adaptive mesh: the code crashes down quickly without BP limiters.  Therefore, BP method and adaptive moving mesh method are both necessary. Figure \ref{Fig:5EQ-RP2} compares the moving mesh and uniform mesh results, where results on moving meshes are obtained by using \[
D\U = \sqrt{\dvert{\rho_x}^2 + \dvert{\rho_{xx}}^2 + 100\dvert{u_x}^2}, 
\]
and $\beta = 0.4$. It is observed from Figure \ref{Fig:5EQ-RP2} that the result produced by the BP moving mesh method has superior resolution than that produced on uniform mesh when $200$ cells are used for the simulation.

\section{Concluding remarks}
In this paper, we developed a simple bound-preserving (BP) moving mesh method for one-dimensional hyperbolic conservation laws. The method first decomposes a  finite-volume scheme into sub-cell schemes, then enforces the BP property by blending high-order sub-cell schemes with their first-order counterparts. Since the BP property depends solely on the first-order sub-cell schemes, it permits significantly larger CFL conditions, thereby substantially improving computational efficiency. We systematically derived the accuracy-preserving (AP) CFL condition that maintains numerical accuracy after BP enforcement. This AP CFL condition is mild, further demonstrating the effectiveness of our approach.

To illustrate the potential of the method, we extended it to the Euler equations, and the five-equation model of two-medium flows involving non-conservative products. Numerical examples validated our theoretical analysis and confirmed the robustness and effectiveness of the method.

Finally, it is worth mentioning that the proposed BP method shows strong promise for high-dimensional problems, a practically important direction for future research. Results will be reported in the forthcoming papers.


\bibliographystyle{siam}
\bibliography{Reference}

@Article{ALLAIRE2002577,
  Title                    = {A Five-Equation Model for the Simulation of Interfaces between Compressible Fluids},
  Author                   = {Grégoire Allaire and Sébastien Clerc and Samuel Kokh},
  Journal                  = {Journal of Computational Physics},
  Year                     = {2002},
  Number                   = {2},
  Pages                    = {577--616},
  Volume                   = {181}
}

@Article{CHENG2014143,
  Title                    = {Positivity-preserving {L}agrangian scheme for multi-material compressible flow},
  Author                   = {Cheng, Juan and Shu, Chi-Wang},
  Journal                  = {Journal of Computational Physics},
  Year                     = {2014},
  Pages                    = {143--168},
  Volume                   = {257},

  ISSN                     = {0021-9991}
}

@Article{CORALIC201495,
  author  = {Coralic, Vedran and Colonius, Tim},
  journal = {Journal of Computational Physics},
  title   = {Finite-volume {WENO} scheme for viscous compressible multicomponent flows},
  year    = {2014},
  pages   = {95-121},
  volume  = {274},
}

@Article{FAN2022111446,
  Title                    = {Positivity-preserving high order finite difference {WENO} schemes for compressible {N}avier-{S}tokes equations},
  Author                   = {Fan, Chuan and Zhang, Xiangxiong and Qiu, Jianxian},
  Journal                  = {Journal of Computational Physics},
  Year                     = {2022},
  Pages                    = {111446},
  Volume                   = {467},

  Doi                      = {https://doi.org/10.1016/j.jcp.2022.111446},
  ISSN                     = {0021-9991},
  Url                      = {https://www.sciencedirect.com/science/article/pii/S0021999122005083}
}

@Article{Gao2023High,
  author  = {Gao, Zhen and Guo, Shuang and Wang, Bao-Shan and Gu, Yaguang},
  journal = {Communications in Computational Physics},
  title   = {High Order Bound- and Positivity-preserving Finite Difference Affine-invariant {AWENO} Scheme for the Five-equation Model of Two-medium Flows},
  year    = {2024},
  number  = {3},
  pages   = {781--820},
  volume  = {36},
}

@Article{GU2023112190,
  author  = {Gu, Yaguang and Gao, Zhen and Hu, Guanghui and Li, Peng and Fu, Qingcheng},
  journal = {Journal of Computational Physics},
  title   = {High Order Well-Balanced Positivity-Preserving Scale-invariant AWENO Scheme for {E}uler Systems with Gravitational Field},
  year    = {2023},
  note    = {Paper No. 112190},
}

@Article{Gu2021A,
  author  = {Gu, Yaguang and Gao, Zhen and Hu, Guanghui and Li, Peng and Wang, Lifeng},
  journal = {Journal of Scientific Computing},
  title   = {A Robust High Order Alternative {WENO} Scheme for the Five-Equation Model},
  year    = {2021},
  volume  = {88},
  notes   = {Paper No. 12},
}

@Article{GU2022An,
  Title                    = {An Adaptive Moving Mesh Method for the Five-Equation Model},
  Author                   = {Gu, Yaguang and Luo, Dongmi and Gao, Zhen and Chen, Yibing},
  Journal                  = {Communications in Computational Physics},
  Year                     = {2022},
  Number                   = {1},
  Pages                    = {189--221},
  Volume                   = {32},

  Doi                      = {https://doi.org/10.4208/cicp.OA-2021-0169},
  ISSN                     = {1991-7120},
  Url                      = {http://global-sci.org/intro/article_detail/cicp/20792.html}
}

@Article{HU2013169,
  Title                    = {Positivity-preserving method for high-order conservative schemes solving compressible {E}uler equations},
  Author                   = {Hu, Xiangyu Y. and Adams, Nikolaus A. and Shu, Chi-Wang},
  Journal                  = {Journal of Computational Physics},
  Year                     = {2013},
  Pages                    = {169--180},
  Volume                   = {242},

  ISSN                     = {0021-9991}
}

@Article{JIANG1996202,
  author  = {Jiang, Guang-Shan and Shu, Chi-Wang},
  journal = {Journal of Computational Physics},
  title   = {Efficient Implementation of Weighted {ENO} Schemes},
  year    = {1996},
  number  = {1},
  pages   = {202-228},
  volume  = {126},
}

@Article{Kurganov2021Adaptive,
  Title                    = {Adaptive Moving Mesh Central-Upwind Schemes for Hyperbolic System of {PDEs}: Applications to Compressible {E}uler Equations and Granular Hydrodynamics},
  Author                   = {Kurganov, Alexander and Qu, Zhuolin and Rozanova, Olga S. and Wu, Tong},
  Journal                  = {Communications on Applied Mathematics and Computation},
  Year                     = {2021},
  Pages                    = {445-479},
  Volume                   = {3},

  Issue                    = {3}
}

@Article{QIAN2018172,
  author  = {Qian, Shouguo and Li, Gang and Shao, Fengjing and Xing, Yulong},
  journal = {Advances in Water Resources},
  title   = {Positivity-preserving well-balanced discontinuous {G}alerkin methods for the shallow water flows in open channels},
  year    = {2018},
  pages   = {172--184},
  volume  = {115},
}

@Article{WANG2022630,
  Title                    = {Affine-invariant {WENO} weights and operator},
  Author                   = {Bao-Shan Wang and Wai Sun Don},
  Journal                  = {Applied Numerical Mathematics},
  Year                     = {2022},
  Pages                    = {630--646},
  Volume                   = {181},

  Doi                      = {https://doi.org/10.1016/j.apnum.2022.07.007},
  ISSN                     = {0168-9274},
  Url                      = {https://www.sciencedirect.com/science/article/pii/S0168927422001817}
}

@Article{WONG2021A,
  author  = {Wong, Man Long and Angel, Jordan B. and Barad, Michael F. and Kiris, Cetin F.},
  journal = {Journal of Computational Physics},
  title   = {A positivity-preserving high-order weighted compact nonlinear scheme for compressible gas-liquid flows},
  year    = {2021},
  note    = {Paper No. 110569},
  volume  = {444},
}

@Article{ZHANG20122245,
  Title                    = {Positivity-preserving high order finite difference {WENO} schemes for compressible {E}uler equations},
  Author                   = {Zhang, Xiangxiong and Shu, Chi-Wang},
  Journal                  = {Journal of Computational Physics},
  Year                     = {2012},
  Number                   = {5},
  Pages                    = {2245-2258},
  Volume                   = {231},

  Doi                      = {https://doi.org/10.1016/j.jcp.2011.11.020},
  ISSN                     = {0021-9991},
  Url                      = {https://www.sciencedirect.com/science/article/pii/S0021999111006759}
}

@Article{ZHANG20111238,
  Title                    = {Positivity-preserving high order discontinuous {G}alerkin schemes for compressible {E}uler equations with source terms},
  Author                   = {Zhang, Xiangxiong and Shu, Chi-Wang},
  Journal                  = {Journal of Computational Physics},
  Year                     = {2011},
  Number                   = {4},
  Pages                    = {1238--1248},
  Volume                   = {230},

  Doi                      = {https://doi.org/10.1016/j.jcp.2010.10.036},
  ISSN                     = {0021-9991},
  Url                      = {https://www.sciencedirect.com/science/article/pii/S0021999110006017}
}

@Article{ZHANG20103091,
  Title                    = {On maximum-principle-satisfying high order schemes for scalar conservation laws},
  Author                   = {Zhang, Xiangxiong and Shu, Chi-Wang},
  Journal                  = {Journal of Computational Physics},
  Year                     = {2010},
  Number                   = {9},
  Pages                    = {3091-3120},
  Volume                   = {229},

  Doi                      = {https://doi.org/10.1016/j.jcp.2009.12.030},
  ISSN                     = {0021-9991},
  Url                      = {https://www.sciencedirect.com/science/article/pii/S0021999109007165}
}

@Article{ZHANG20108918,
  Title                    = {On positivity-preserving high order discontinuous {G}alerkin schemes for compressible {E}uler equations on rectangular meshes},
  Author                   = {Zhang, Xiangxiong and Shu, Chi-Wang},
  Journal                  = {Journal of Computational Physics},
  Year                     = {2010},
  Number                   = {23},
  Pages                    = {8918--8934},
  Volume                   = {229},

  Doi                      = {https://doi.org/10.1016/j.jcp.2010.08.016},
  ISSN                     = {0021-9991},
  Url                      = {https://www.sciencedirect.com/science/article/pii/S0021999110004535}
}

@Article{FU2022111600,
  author  = {Fu, Pei and Xia, Yinhua},
  journal = {Journal of Computational Physics},
  title   = {The positivity preserving property on the high order arbitrary {L}agrangian-{E}ulerian discontinuous {G}alerkin method for {E}uler equations},
  year    = {2022},
  note    = {Paper No. 111600},
}

@Article{Tang2003Adaptive,
  author  = {Tang, Huazhong and Tang, Tao},
  journal = {SIAM Journal on Numerical Analysis},
  title   = {Adaptive Mesh Methods for One- and Two-Dimensional Hyperbolic Conservation Laws},
  year    = {2003},
  number  = {2},
  pages   = {487-515},
  volume  = {41},
}

@Article{CUI2025114189,
  author  = {Cui, Shumo and Gu, Yaguang and Kurganov, Alexander and Xin, Ruixiao and Wu, Kailiang},
  journal = {Journal of Computational Physics},
  title   = {Positivity-preserving new low-dissipation central-upwind schemes for compressible {E}uler equations},
  year    = {2025},
  note    = {Paper No. 114189},
  volume  = {538},
}

@Article{Li2006Moving,
  author  = {Li, Ruo and Tang, Tao},
  journal = {Journal of Scientific Computing},
  title   = {Moving Mesh Discontinuous {G}alerkin Method for Hyperbolic Conservation Laws},
  year    = {2006},
  pages   = {347 - 363},
  volume  = {27},
  issue   = {1},
}

@Article{Gu2017A,
  author  = {Gu, Yaguang and Hu, Guanghui},
  journal = {Communications in Computational Physics},
  title   = {A Third Order Adaptive {ADER} Scheme for One Dimensional Conservation Laws},
  year    = {2017},
  issn    = {1991-7120},
  number  = {3},
  pages   = {829--851},
  volume  = {22},
  doi     = {https://doi.org/10.4208/cicp.OA-2016-0088},
  url     = {http://global-sci.org/intro/article_detail/cicp/9983.html},
}

@Article{Kurganov2022Well,
  author  = {Kurganov, Alexander and Qu, Zhuolin and Wu, Tong},
  journal = {ESAIM: M2AN},
  title   = {Well-balanced positivity preserving adaptive moving mesh central-upwind schemes for the {S}aint-{V}enant system},
  year    = {2022},
  number  = {4},
  pages   = {1327-1360},
  volume  = {56},
  doi     = {10.1051/m2an/2022041},
  url     = {https://doi.org/10.1051/m2an/2022041},
}

@Article{PAN2020109558,
  author   = {Pan, Liang and Zhao, Fengxiang and Xu, Kun},
  journal  = {Journal of Computational Physics},
  title    = {High-order ALE gas-kinetic scheme with {WENO} reconstruction},
  year     = {2020},
  issn     = {0021-9991},
  pages    = {109558},
  volume   = {417},
  doi      = {https://doi.org/10.1016/j.jcp.2020.109558},
  url      = {https://www.sciencedirect.com/science/article/pii/S0021999120303326},
}

@Article{Fu2019Arbitrary,
  author  = {Fu, Pei and Schn\"ucke, Gero and Xia, Yinhua},
  journal = {Mathematics of Computation},
  title   = {Arbitrary {L}agrangian-{E}ulerian discontinuous {G}alerkin method for conservation laws on moving simplex meshes},
  year    = {2019},
  issn    = {0025-5718},
  number  = {319},
  pages   = {2221--2255},
  volume  = {88},
}

@Book{Kuzmin2023Property,
  author    = {Kuzmin, Dmitri and Hajduk, Hennes},
  publisher = {WORLD SCIENTIFIC},
  title     = {Property-preserving numerical schemes for conservation laws},
  year      = {2023},
}

@Article{Fu2025Bound,
  author  = {Fu, Qingcheng and Gu, Yaguang and Kurganov, Alexander and Wang, Bao-Shan},
  journal = {Journal of Scientific Computing},
  title   = {Bound- and Positivity-Preserving Path-Conservative Central-Upwind {AWENO} Scheme for the Five-Equation Model of Compressible Two-Component Flows},
  year    = {2025},
  note    = {Paper No. 94},
  volume  = {104},
}

@Article{KURGANOV2023111773,
  author  = {Kurganov, Alexander and Liu, Yongle and Xin, Ruixiao},
  journal = {Journal of Computational Physics},
  title   = {Well-balanced path-conservative central-upwind schemes based on flux globalization},
  year    = {2023},
  volume  = {474},
  notes   = {Paper No. 111773},
}

@Article{XU2014Parametrized,
  author  = {Xu, Zhengfu},
  journal = {Mathematics of Computation},
  title   = {Parametrized maximum principle preserving flux limiters for high order schemes solving hyperbolic conservation laws: one-dimensional scalar problem},
  year    = {2014},
  pages   = {2213-2238},
  volume  = {83},
}

@Article{Wu2021Uniformly,
  author  = {Wu, Kailiang and Xing, Yulong},
  journal = {SIAM Journal on Scientific Computing},
  title   = {Uniformly High-Order Structure-Preserving Discontinuous {G}alerkin Methods for {E}uler Equations with Gravitation: Positivity and Well-Balancedness},
  year    = {2021},
  number  = {1},
  pages   = {A472-A510},
  volume  = {43},
}

@Article{BORIS197338,
  author  = {Boris, Jay P and Book, David L},
  journal = {Journal of Computational Physics},
  title   = {Flux-corrected transport. I. {SHASTA}, a fluid transport algorithm that works},
  year    = {1973},
  number  = {1},
  pages   = {38-69},
  volume  = {11},
}

@Book{kuzmin2012flux,
  author    = {Kuzmin, Dmitri and L{\"o}hner, Rainald and Turek, Stefan},
  publisher = {Springer Science \& Business Media},
  title     = {Flux-corrected transport: principles, algorithms, and applications},
  year      = {2012},
}

@Article{Xiong2016Parametrized,
  author  = {Xiong, Tao and Qiu, Jing-Mei and Xu, Zhengfu},
  journal = {Journal of Scientific Computing},
  title   = {Parametrized Positivity Preserving Flux Limiters for the High Order Finite Difference {WENO} Scheme Solving Compressible {E}uler Equations},
  year    = {2016},
  pages   = {1066-1088},
  volume  = {67},
  issue   = {3},
}

@Article{Liang2014Parametrized,
  author  = {Liang, Chao and Xu, Zhengfu},
  journal = {Journal of Scientific Computing},
  title   = {Parametrized Maximum Principle Preserving Flux Limiters for High Order Schemes Solving Multi-Dimensional Scalar Hyperbolic Conservation Laws},
  year    = {1996},
  pages   = {41-60},
  volume  = {58},
}

@Article{KUZMIN2020112804,
  author  = {Kuzmin, Dmitri},
  journal = {Computer Methods in Applied Mechanics and Engineering},
  title   = {Monolithic convex limiting for continuous finite element discretizations of hyperbolic conservation laws},
  year    = {2020},
  note    = {Paper No. 112804},
  volume  = {361},
}

@Article{BOOK1975248,
  author  = {Book, David and Boris, Joao and Hain, K.},
  journal = {Journal of Computational Physics},
  title   = {Flux-corrected transport {II}: Generalizations of the method},
  year    = {1975},
  number  = {3},
  pages   = {248-283},
  volume  = {18},
}

@Article{BORIS1976397,
  author  = {Joao Boris and David Book},
  journal = {Journal of Computational Physics},
  title   = {Flux-corrected transport. {III}. Minimal-error FCT algorithms},
  year    = {1976},
  number  = {4},
  pages   = {397-431},
  volume  = {20},
}

@Article{ZALESAK1979335,
  author  = {Zalesak, Steven},
  journal = {Journal of Computational Physics},
  title   = {Fully multidimensional flux-corrected transport algorithms for fluids},
  year    = {1979},
  number  = {3},
  pages   = {335-362},
  volume  = {31},
}

@Article{Abgrall2025bound,
  author  = {Abgrall, R\'{e}mi and Jiao, Miaosen and Liu, Yongle and Wu, Kailiang},
  journal = {Communications in Computational Physics},
  title   = {Bound-preserving {P}oint-{A}verage-{M}oment {P}olynomi{A}l-interpreted ({PAMPA}) scheme: one-dimensional case},
  year    = {2025},
  number  = {1},
  pages   = {29-58},
  volume  = {39},
}

@Article{Duan2025Active,
  author  = {Duan, Junming and Barsukow, Wasilij and Klingenberg, Christian},
  journal = {SIAM Journal on Scientific Computing},
  title   = {Active Flux Methods for Hyperbolic Conservation Laws-Flux Vector Splitting and Bound-Preservation},
  year    = {2025},
  number  = {2},
  pages   = {A811-A837},
  volume  = {47},
}

@Article{Guermond2017Invariant,
  author  = {Guermond, Jean-Luc and Popov, Bojan},
  journal = {SIAM Journal on Numerical Analysis},
  title   = {Invariant Domains and Second-Order Continuous Finite Element Approximation for Scalar Conservation Equations},
  year    = {2017},
  number  = {6},
  pages   = {3120-3146},
  volume  = {55},
}

@Article{GOTTLIEB1998Total,
  author  = {Gottlieb, Sigal AND Shu, Chi-Wang},
  journal = {Mathematics of Computation},
  title   = {Total Variation Diminishing {R}unge-{K}utta Schemes},
  year    = {1998},
  number  = {221},
  pages   = {73-85},
  volume  = {67},
}

@Article{Shu2001Strong,
  author  = {Shu, Chi-Wang},
  journal = {SIAM Review},
  title   = {Strong Stability-Preserving High-Order Time Discretization Methods},
  year    = {2001},
  pages   = {89-112},
  volume  = {43},
}

@InProceedings{Tang2005Moving,
  author    = {Tang, Tao},
  booktitle = {Recent Advances in Adaptive Computation},
  title     = {Moving mesh methods for computational fluid dynamics},
  year      = {2005},
  editor    = {Shi, Z.-C. and Chen, Z. and Tang, T. and Yu, D.},
  note      = {Proceedings of the International Conference on Recent Advances in Adaptive Computation, May 2004, Hangzhou, China},
  pages     = {141-173},
  publisher = {American Mathematical Society},
  volume    = {383},
}

\end{document}